\documentclass{amsart}

\usepackage[dvipsnames]{xcolor}

\usepackage{amsthm,amssymb,amsmath,amsfonts,mathrsfs,amscd,yfonts,stmaryrd}
\usepackage{url,tikz,enumitem,caption,stmaryrd,comment,mathrsfs,marginnote}
\usepackage[all]{xy}
\usepackage{latexsym}
\usepackage[utf8]{inputenc}
\usepackage[T1]{fontenc}
\usepackage{tikz-cd}
\usepackage{libertine}
\usetikzlibrary{matrix,positioning,arrows,shapes,chains,calc,automata}
\usetikzlibrary{decorations.pathmorphing}
\usetikzlibrary{arrows,automata,positioning}
\usepackage[pagebackref=true]{hyperref}

\newlength{\proofmargin}
\renewcommand*{\backref}[1]{}
\renewcommand*{\backrefalt}[4]{%
  \ifcase #1 %
    \relax
  \or
    $\uparrow$#2.%
  \else
    $\uparrow$#2.%
  \fi%
}

  \renewenvironment{thebibliography}[1]{%
   \begin{oldthebibliography}{#1}%
      \setlength{\parskip}{0ex}%
      \setlength{\itemsep}{0.1ex}%
       \setlength{\labelwidth}{1.4cm} 
  }%
  {%
    \end{oldthebibliography}%
  }

\newcommand{\C}{\mathbb C}
\newcommand{\R}{\mathbb R}

\newcommand{\res}{\mathrm{res}}
\newcommand{\Z}{\mathbb Z}
\newcommand{\Q}{\mathbb Q}

\newcommand{\PP}{\mathbb P}
\newcommand{\JJ}{\mathcal J}

\newcommand{\cH}{\mathcal H}
\newcommand{\cO}{\mathcal O}
\newcommand{\cF}{\mathcal F}
\newcommand{\cE}{\mathcal E}
\newcommand{\cC}{\mathcal C}
\newcommand{\cD}{\mathcal D}
\newcommand{\cT}{\mathcal T}
\newcommand{\cR}{\mathcal R}

\newcommand{\fm}{\mathfrak m}
\newcommand{\fl}{\mathfrak l}
\newcommand{\fd}{\mathfrak d}
\newcommand{\fa}{\mathfrak a}
\newcommand{\fp}{\mathfrak p}

\newcommand{\red}{\mathrm{red}}
\newcommand{\loc}{\mathrm{loc}}

\newcommand{\Tr}{\mathrm{Tr}}
\newcommand{\ad}{\mathrm{ad}}

\newcommand{\eps}{\varepsilon}
\newcommand{\ab}{\mathrm{ab}}
\newcommand{\rec}{\mathrm{rec}}

\newcommand{\Hom}{\mathrm{Hom}}
\newcommand{\Gal}{\mathrm{Gal}}
\newcommand{\GL}{\mathrm{GL}}

\newcommand{\rH}{\mathrm H}

\newcommand{\ord}{\mathrm{ord}}
\newcommand{\nord}{\mathrm{n.ord}}

\newcommand{\Spm}{\mathrm{Spm}}

\newcommand{\Frob}{\mathrm{Frob}}

\newcommand{\rZ}{\mathrm Z}

\newcommand{\univ}{\mathrm{univ}}
\newcommand{\sL}{\mathscr {L}}
\newcommand{\sbL}{\mathscrbf {L}}
\renewcommand{\phi}{\varphi}
\newcommand{\TT}{{\mathbb T}}
\newcommand{\fN}{{\mathfrak N}}
\newcommand{\bfk}{{\mathbf k}}
\newcommand{\bfw}{{\mathbf w}}
\newcommand{\bfv}{{\mathbf v}}
\newcommand{\bft}{{\mathbf t}}

\hypersetup{linktocpage,bookmarksopen = true, bookmarksopenlevel = 1, colorlinks=true,
linkcolor = Cerulean, citecolor = RubineRed , urlcolor = gray, pdfstartview = FitR} 

\DeclareMathOperator{\SL}{\mathrm{SL}}
\DeclareMathOperator{\cM}{\mathcal{M}}
\DeclareMathOperator{\cA}{\mathcal{A}}
\DeclareMathOperator{\dlog}{\mathrm{dlog}}
\DeclareMathOperator{\lra}{\longrightarrow}
\DeclareMathOperator{\Cl}{\mathrm{Cl}}

\DeclareMathOperator{\bX}{\mathbb{X}}
\DeclareMathOperator{\bD}{\mathbb{D}}
\DeclareMathOperator{\LC}{\mathrm{LC}}
\DeclareMathOperator{\Nm}{\mathrm{Nm}}
\DeclareMathOperator{\A}{\mathbf{A}}
\DeclareMathOperator{\ft}{\mathfrak{t}}

\DeclareMathOperator{\Stab}{\mathrm{Stab}}
\DeclareMathOperator{\RM}{\mathrm{RM}}
\DeclareMathOperator{\bfU}{{\mathsf U}}
\DeclareMathOperator{\bfG}{{\mathsf G}}
\DeclareMathOperator{\bfZ}{{\mathsf Z}}
\DeclareMathOperator{\bfUtor}{{\mathsf U}^{\circ}}

\DeclareMathOperator{\bfZtor}{{\mathsf Z}^{\circ}}
\DeclareMathAlphabet{\mathbbm}{U}{bbm}{m}{n}
\DeclareMathAlphabet{\mathscrbf}{OMS}{mdugm}{b}{n} 

\newcommand{\mat}[4]{\left(\begin{array}{cc}#1&#2\\#3&#4\end{array}\right)}
\newcommand{\smallmat}[4]{\bigl(\begin{smallmatrix}#1&#2\\#3&#4\end{smallmatrix}\bigr)}
\newcommand{\mint}{\times\!\!\!\!\!\!\!\hspace{.097em}\int}
\newcommand{\PGL}{{\mathrm {PGL}}}

\makeatletter
\DeclareFontFamily{OMX}{MnSymbolE}{}
\DeclareSymbolFont{MnLargeSymbols}{OMX}{MnSymbolE}{m}{n}
\SetSymbolFont{MnLargeSymbols}{bold}{OMX}{MnSymbolE}{b}{n}
\DeclareFontShape{OMX}{MnSymbolE}{m}{n}{
    <-6>  MnSymbolE5
   <6-7>  MnSymbolE6
   <7-8>  MnSymbolE7
   <8-9>  MnSymbolE8
   <9-10> MnSymbolE9
  <10-12> MnSymbolE10
  <12->   MnSymbolE12
}{}
\DeclareFontShape{OMX}{MnSymbolE}{b}{n}{
    <-6>  MnSymbolE-Bold5
   <6-7>  MnSymbolE-Bold6
   <7-8>  MnSymbolE-Bold7
   <8-9>  MnSymbolE-Bold8
   <9-10> MnSymbolE-Bold9
  <10-12> MnSymbolE-Bold10
  <12->   MnSymbolE-Bold12
}{}

\let\llangle\@undefined
\let\rrangle\@undefined
\let\lsem\@undefined
\let\rsem\@undefined
\DeclareMathDelimiter{\llangle}{\mathopen}%
                     {MnLargeSymbols}{'164}{MnLargeSymbols}{'164}
\DeclareMathDelimiter{\rrangle}{\mathclose}%
                     {MnLargeSymbols}{'171}{MnLargeSymbols}{'171}
\DeclareMathDelimiter{\lsem}{\mathopen}%
                     {MnLargeSymbols}{'102}{MnLargeSymbols}{'102}
\DeclareMathDelimiter{\rsem}{\mathclose}%
                     {MnLargeSymbols}{'107}{MnLargeSymbols}{'107}
\makeatother

\newtheorem{theorem}{Theorem}[section]
\newtheorem{lemma}[theorem]{Lemma}

\newtheorem{proposition}[theorem]{Proposition}
\newtheorem{definition}[theorem]{Definition}

\newtheorem{prop}[theorem]{Proposition}

\theoremstyle{remark}

\newtheorem{remarkwr}[table]{Remark}
\newenvironment{remark}{\begin{remarkwr}\begin{upshape}}{\end{upshape}\end{remarkwr}}

\makeatletter
\renewenvironment{proof}[1][\proofname]%
{%
\par\pushQED{\qed}\normalfont\topsep6\p@\@plus6\p@\relax%
\begin{list}{}{\rightmargin=8pt\leftmargin=\proofmargin}%
  \item[\hskip\labelsep\bfseries#1\@addpunct{.}]\ignorespaces
}{%
\popQED\end{list}\@endpefalse%
}%
\makeatother

\title[The RM values of the Dedekind--Rademacher cocycle]{The  values of  the Dedekind--Rademacher cocycle   \\ at real multiplication points}
\author{Henri Darmon, Alice Pozzi and Jan Vonk}
\date{}
\address{}
\email{}

\include{thebibliography}

\address{H. D.: Montreal, Canada}
\email{darmon@math.mcgill.ca}
\address{A. P.: London, UK}
\email{alice.pozzi89@gmail.com}
\address{J.V.: Leiden, Netherlands}
\email{j.b.vonk@math.leidenuniv.nl}


\begin{document}

\begin{abstract}
The values of the so-called {\em Dedekind--Rademacher cocycle} at certain real quadratic arguments are 
 shown to  be    global $p$-units  in the narrow Hilbert 
class field  of the associated  real quadratic field, as predicted by the conjectures of \cite{darmon-dasgupta} and \cite{darmon-vonk}.  
The strategy for proving  this result  combines the approach 
of \cite{DPV1} with one   crucial extra ingredient:  
 the study of infinitesimal deformations of irregular 
 Hilbert Eistenstein series of weight one in the anti-parallel
  direction, building on the techniques  of \cite{betina-dimitrov-pozzi}.  
\end{abstract}

\maketitle
\tableofcontents

\section*{Introduction}

Let $\cH_p$ denote Drinfeld's $p$-adic upper half plane, and let $ \cM^\times$ 
denote the multiplicative group  of non-zero rigid meromorphic functions on 
$\cH_p$, equipped with the translation action of the discrete  group
$\SL_2(\Z[1/p])$ by M\"obius transformations. A {\em rigid meromorphic cocycle} on a congruence subgroup $\Gamma\subset \SL_2(\Z[1/p])$ is a class in $\rH^1(\Gamma, \cM^\times)$. If $\tau\in \cH_p$ is a {\em real multiplication}, or RM, point, i.e.,  generates a real quadratic extension of $\Q$, the {\em value} of $J$ at $\tau$ is defined to be
\begin{equation}
\label{eqn:rm-value}
 J[\tau] := J(\gamma_\tau)(\tau) \in \C_p\cup \{\infty\},
 \end{equation}
 where $\gamma_\tau \in \Gamma$ is the {\em automorph} of $\tau$,
  a suitably normalised generator of the stabiliser of $\tau$ in 
 $\Gamma$.
The relevance of the  RM values of rigid meromorphic cocycles    to explicit class field theory for real quadratic fields 
 has been explored in \cite{darmon-hpxh},
\cite{darmon-dasgupta}  \cite{darmon-vonk}, and 
\cite{darmon-vonk-borcherds},
where  it is
 conjectured, broadly speaking, that they behave in many key respects just  like the values of classical modular functions at CM points, and in particular that they belong to, and often generate, narrow ring class fields  of real quadratic fields.

\par Theorem B below gives some theoretical evidence for this general conjecture in the simplest case where $\Gamma = \SL_2(\Z[1/p])$ and $J$ is {\em  analytic}, i.e., takes values in the subgroup $\cA^\times\subset \cM^\times$ of rigid analytic functions. Strictly speaking, there are no interesting  rigid analytic cocycles: the group $\rH^1(\Gamma, \cA^\times)$ is 
generated, up to torsion, by the class $J_{{\rm triv}}$ given by
$$ J_{{\rm triv}}\mat{a}{b}{c}{d}(z) = cz+d,$$
whose RM values are units in the associated quadratic order -- hence, algebraic, but not in an interesting way for explicit class field theory. 

\par There is a  less trivial class in $\rH^1(\Gamma, \cA^{\times}/p^{\Z})$ arising from the classical {\em Dedekind--Rademacher homomorphism} $\varphi_{\rm DR}: \Gamma_0(p) \lra \Z$ 
 describing the periods of the weight two Eisenstein series  
     \begin{equation}
    \label{eqn:E2p}
    E_2^{(p)}(q) = \dlog\left(\frac{\Delta(q^p)}{\Delta(q)}\right) =
    \left(  p-1  +  24 \sum_{n=1}^\infty \sigma^{(p)}(n) q^n\right)\frac{dq}{q}, 
    \qquad  \mbox{ where } \sigma^{(p)}(n) := \sum_{p\nmid d|n} d,
    \end{equation}
and given by
\begin{equation}
\label{eqn:varphi-DR}
\varphi_{\rm DR}(\gamma) := \frac{1}{2\pi i} \int_{z_0}^{\gamma z_0}   \!\!  2 E_2^{(p)}(z) dz.
\end{equation}
More precisely, the description of $\Gamma$ as an amalgamated product of two conjugate copies of $\SL_2(\Z)$ intersecting in $\Gamma_0(p)$ leads to  an injection 
\[
\rH^1(\Gamma_0(p),\Z) \hookrightarrow \rH^2(\Gamma,\Z).
\]
Let $\alpha_{_{\rm DR}} \in Z^2(\Gamma,\Z)$ be a two-cocycle whose cohomology class is the image of $\varphi_{\rm DR}$ under this map. Refining a construction
 of \cite{darmon-dasgupta},  Theorem A below asserts  that the cocycle $p^{\alpha_{\rm DR}}$  with values in $p^\Z$ is trivialised in the larger  group   $\cA^\times \supset p^\Z$:

\medskip\medskip
\noindent
{\bf Theorem A}. {\em 
There is a one-cochain $J_{\rm DR}\in C^1(\Gamma, \cA^\times)$ satisfying
\[
\gamma_1 J_{\rm DR}(\gamma_2) \div J_{\rm DR}(\gamma_1\gamma_2) \times J_{\rm DR}(\gamma_1) = p^{\alpha_{_{\rm DR}}(\gamma_1,\gamma_2)}, 
\quad \mbox{ for all } \gamma_1, \gamma_2 \in \Gamma.
\]  
}

\medskip
The essential triviality of $\rH^1(\Gamma,\cA^\times)$ shows that $J_{\rm DR}$ is uniquely determined  up to coboundaries and powers of the cocycle $J_{\rm triv}$ above. The proof of Theorem A is  given in \S~\ref{sec:dedekind-rademacher}, and constructs an explicit cochain $J_{\rm DR}$ which is well defined up to coboundaries. The natural image of $J_{\rm DR}$ in $\rH^1(\Gamma,\cA^\times/p^{\Z})$ is the {\em Dedekind--Rademacher cocycle} of the title. The rigid analytic cocycles  of higher level studied in  \cite{darmon-dasgupta} are all multiplicative combinations of $\GL_2(\Q)$-translates of this basic cocycle. The  proof of Theorem A  complements the approach of \cite{darmon-dasgupta}, producing a more canonical object in level $1$ which can be envisaged as an avatar of the Eisenstein series $E_2$ in the setting of rigid meromorphic cocycles. The RM values of  $J_{\rm DR}$ are well defined modulo $p^{\Z}$, and it therefore makes sense to enquire about their algebraicity, and their factorisation away from $p$. 
  
\par An RM point $\tau\in \cH_p$ is said to be of {\em discriminant } $D$ if it satisfies an equation of the form $Q(\tau,1) = 0$, where $Q(x,y) = Ax^2 + B xy + Cy^2$ is a primitive integral binary quadratic form of discriminant $D$. The set $\cH_p^D$ of $\tau$ of  a fixed discriminant $D$ is non-empty precisely when $p$ is inert or ramified in the quadratic field $F = \Q(\sqrt{D})$, and   is preserved by the action of $\SL_2(\Z)$. The orbit set $\SL_2(\Z) \backslash \cH_p^D$ is in natural bijection with  the class group $\Cl(D)$, by sending the orbit of $\tau$ to the narrow equivalence class $\cC_\tau$ of the fractional ideal generated by $\tau$ and $1$ when $\tau -\tau'$ is positive. The reciprocity map 
\[
\rec:\Cl(D) \lra \Gal(H/F)
\]
of global class field theory identifies $\Cl(D)$ with the Galois group of the narrow ring  class field of $H$ over $F$ attached to $D$. If $\cH_p^D$ is non-empty and $p\nmid D$, then the prime $p$ is inert in $F/\Q$ and  splits completely in $H/F$. The choice of an embedding $\bar\Q\subset \bar\Q_p$ hence determines a prime $\fp$ of $H$ above $p$, which is fixed once and for all. Fix also a complex embedding $\bar\Q\subset \C$ and write $x\mapsto \bar x$ for the action of complex conjugation on $H$ (which is independent of the choice of embedding). Let 
\[
\cO_H[1/p]^\times_-
\] 
be the group of $p$-units of $H$ which are in the minus-eigenspace for the action of complex conjugation. By the Dirichlet $S$-unit theorem, it is a $\Z$-module of rank $[H:F]/2$ if 
$F$  does not possess a unit of negative norm, and  is finite otherwise. In particular, there is a unique element $u_\tau\in (\cO_H[1/p]^\times_-)\otimes \Q$ satisfying 
 \begin{equation}
\label{eqn:factorisation-JDR-tau}
{\rm ord}_{\fp^\sigma}(u_\tau) = - L(F,\cC_\tau^\sigma,0), \quad
\mbox{ for all } \sigma\in {\rm Gal}(H/F),
\end{equation}
where
$L(F,\cC_\tau^\sigma,s)$ is the partial zeta function of the narrow  ideal class $\cC_\tau^\sigma$ (cf.~\cite[Prop.~3.8]{gross-padic}). The $p$-unit $u_\tau$ is called the {\em Gross--Stark unit} attached to $H/F$ (and the prime $\fp$). The Brumer--Stark conjecture implies that $u_\tau^{12}$  belongs to $\cO_H[1/p]^\times$ rather than to  the tensor product of this group with $\Q$. The proof by  Samit Dasgupta and Mahesh Kakde of  (the  prime to $2$ part of) the  Brumer--Stark conjecture \cite{DK-brumer}   in this setting shows that $u_\tau^{12}$ belongs to $(\cO_H[1/p]^\times)\otimes \Z[1/2]$. 
 
\par The principal conjecture of \cite{darmon-dasgupta}, and its  refinement covering the Dedekind--Rademacher cocycle itself, asserts that $J_{\rm DR}[\tau]$ is equal, 
up to a small torsion ambiguity and powers of $p$,
 to  an integer power of the Gross--Stark unit $u_\tau$.
 The weaker equality 
  \begin{equation}
 \label{eqn:ddp}
  {\rm Norm}_{\Q_{p^2}/\Q_p} (J_{\rm DR}[\tau]) = {\rm Norm}_{\Q_{p^2}/\Q_p}(u_\tau^{12})
 \pmod{(\Q_{p}^{\times})_{\rm tors}, p^{\Z}}
 \end{equation}
 involving the norms to $\Q_p^\times$ of these invariants
 was shown in \cite{darmon-dasgupta} to follow from Gross's $p$-adic analogue of the Stark conjecture on $p$-adic abelian $L$-series of totally real fields at $s=0$  -- at least, after replacing $J_{\rm DR}[\tau]$ by the closely allied quantities denoted $u(\alpha,\tau)$ in 
  \cite{darmon-dasgupta}, which depend on the choice of a suitable modular unit $\alpha \in \cO_{Y_1(N)}^\times$ with auxiliary  level structure.
 The Gross--Stark
 conjecture was then proved in \cite{DDP11}. An important recent work of Samit Dasgupta and Mahesh Kakde  \cite{dasgupta-kakde} has significantly refined the approach  of \cite{DDP11} 
 to prove Gross's {\em tame refinement } of the Gross--Stark conjecture, for arbitrary totally real fields.
Specialising this result to the case of a real  quadratic field 
 leads   to  the refinement 
\begin{equation}
\label{eqn:dk}
J_{\rm DR}[\tau] = u_\tau^{12} \pmod{  (\Q_{p^2}^{\times})_{\rm tors}, \ \ p^\Z}
\end{equation}
of \eqref{eqn:ddp} in which the norm is removed. The removal of this ambiguity  is 
crucial for a truly satisfying approach
 to  explicit class field theory for real quadratic fields.

\par The  main contribution of this  paper is an independent and more direct proof of  \eqref{eqn:dk} for fundamental discriminants: 

\medskip
\noindent
 {\bf Theorem B}. {\em 
Let $D>0$ be a fundamental
 discriminant that is prime to $p$.
 If $\tau$ is an RM point in $\cH_p$ of discriminant $D$, 
 then  $J_{\rm DR}[\tau]$  is equal to the Gross--Stark unit $u_\tau^{12}$, up to torsion in $\Q_{p^2}$ and powers of
 $p$, and in particular belongs to $(\cO_H[1/p]^{\times}) \otimes \Z[1/2]$.  }

\medskip\medskip 
To situate the approach of this paper in the context of previous works, note that Dasgupta and Kakde tackle Theorem B by studying Mazur--Tate style ``tame  refinements" of  the techniques of \cite{DDP11}, leading to a proof of Gross's  tame refinement of his $p$-adic Stark conjecture (known as the ``tower of fields conjecture" \cite{gross-tame}). They then show that this tame refinement implies Theorem B.  Like \cite{dasgupta-kakde}, the present work rests on the 
 careful study of   deformations of Galois representations that was also exploited in \cite{DDP11}, but otherwise differs in its approach to Theorem B by avoiding the recourse to tame deformations.  Its key idea is to package the RM values of $J_{\rm DR}$ as the coefficients of certain modular generating series. The resulting identities (cf.~Theorem C below) are of interest in their own right  and enrich the tapestry of analogies between RM values of rigid meromorphic cocycles and  
CM values of modular functions.

\par The Dedekind--Rademacher cocycle, taken modulo $\C_p^\times$ rather than $p^{\Z}$, is a prototypical instance of a {\em rigid analytic  theta-cocycle}: a function $J:\Gamma \lra \cA^\times$ which satisfies the one-cocycle relation, but only up to multiplicative scalars. The proof of Theorem B rests on the study of  another theta-cocycle, the so-called {\em winding cocycle}  
\begin{equation}
 J_w  \in \rH^1(\Gamma,\cA^{\times}/\C_p^{\times}),
\end{equation}
whose key properties are recalled in \S~\ref{sec:winding}. The notion of RM value can be extended to  theta cocycles by noting that, if the RM point $\tau$ has discriminant prime to $p$, then its automorph $\gamma_\tau$ belongs to $\SL_2(\Z)$. The groups $\rH^1(\SL_2(\Z),\C_p^{\times})$ and $\rH^2(\SL_2(\Z),\C_p^{\times})$ are finite of order dividing $12$, which implies that the restriction of $J^{12}$ to $\SL_2(\Z)$ admits an essentially unique lift $\tilde J\in \rH^1(\SL_2(\Z),\cA^\times)$, and the value $J[\tau]$ can then be defined as in \eqref{eqn:rm-value}, with $J$ replaced by $\tilde J^{1/12}$ on the right hand side. Although there is some torsion ambiguity in the resulting RM values, the $p$-adic logarithms of these RM values are well-defined.

\par The explicit nature of $J_w$ can be parlayed into a proof of the following result:
 
\medskip
\noindent
{\bf Theorem C}. {\em 
Let $\tau$ be as in Theorem B. There is a classical modular form $G_\tau$ of weight two on $\Gamma_0(p)$ with $p$-adic Fourier coefficients, whose $q$-expansion is given by
\[
G_\tau(q) =  \log(u_\tau) + \sum_{n=1}^\infty \log((T_n J_w)[\tau]) q^n,
\]
where $\log \ : \ \cO_{\C_p}^\times \lra \C_p$ is the $p$-adic logarithm. The modular form $G_\tau$ is non-trivial if and only if $\Q(\sqrt{D})$ does not admit a unit of norm $-1$.
}

\par The modular generating series of Theorem C is constructed from the diagonal restriction of a nearly ordinary  deformation of a weight one Hilbert Eisenstein series for $\SL_2(\cO_F)$ in the anti-parallel direction. The logarithm of the global $p$-unit  $u_\tau$ enters into the proof as the eigenvalue of the Frobenius at $p$ on a quotient of the associated $p$-adic Galois representation, via a calculation which exploits the reciprocity law of global class field theory, thereby leveraging  class field theory for $H$ into explicit class field theory for $F$. An essential ingredient in the proof of Theorem C is the study of 
$p$-adic deformations of  irregular Hilbert Eisenstein series  of weight one, which is explained in \S~\ref{sec:deformation} and  forms the technical core of this article.  
This approach is inspired by the study  of 
the local geometry of the  modular eigenvariety in the neighbourhood of irregular Eisenstein points of weight one carried out  in 
  \cite{betina-dimitrov-pozzi}, and its extension 
  to the Hilbert setting   in
 \cite{BDS}.

\par Derivatives of $p$-adic families of (classical, or Hilbert) modular forms can be viewed as  $p$-adic counterparts of incoherent Eisenstein series in the sense of Kudla, and provide a protoypical instance of what might be envisaged as {\em $p$-adic mock modular forms}. Deformations of weight one Hilbert modular Eisenstein series in the parallel weight direction and their diagonal restrictions are studied in \cite{DPV1}, where they are related to the {\em norms to $\Q_p$}  of  $J_{\rm DR}[\tau]$. Because of the loss of information inherent in taking  the norm, Theorem C  represents a significant strengthening of the main theorem of \cite{DPV1}, just as Theorem B strengthens the equality \eqref{eqn:ddp} resulting from the proof of Gross--Stark conjecture in the setting of odd ring class characters of  real quadratic fields. 

\par In \S~\ref{sec:dedekind-rademacher} the Dedekind--Rademacher cocycle is constructed, thereby proving Theorem A. The definition and main properties of the winding cocycle appear in \S~\ref{sec:winding}, where Theorem B is reduced to Theorem C. The modular generating series $G_\tau$ of Theorem C is constructed in \S~\ref{sec:deformation}--\ref{sec:generating-series}. The pivotal \S~\ref{sec:deformation} studies infinitesimal $p$-adic deformations of weight one Hilbert Eisenstein series and their Fourier expansions. Finally, through a calculation carried out in \S~\ref{sec:generating-series}, the form $G_\tau$ is obtained from the ordinary projection of the diagonal restriction of this infinitesimal deformation.

\section{The Dedekind--Rademacher cocycle}
\label{sec:dedekind-rademacher}

This section constructs a one-cochain satisfying Theorem A, which is well-defined up to coboundaries and whose image in $\rH^1(\Gamma,\cA^\times/p^{\Z})$ is the Dedekind--Rademacher cocycle $J_{\rm DR}$ of the introduction.

\subsection{Siegel units}
Let $\cO_{\cH}^\times$ denote the multiplicative group of nowhere vanishing holomorphic functions on the Poincar\'e 
upper half-plane, endowed with the  right  action of $\SL_2(\R)$ given by
$$h|\gamma(z) = h(\gamma z),$$
where $\gamma z$ denotes the usual action of $\gamma$ by M\"obius transformations.

The construction of $J_{\rm DR}$ rests  on the   Siegel units ${_c}g_{\alpha,\beta} \in \cO_{\cH}^\times$ 
indexed by pairs $ (\alpha,\beta)\in (\Q/\Z)^2 - \{(0,0)\}$ of order $N>1$, 
depending on an auxiliary integer 
$c$ which is relatively prime to $6N$. 
They satisfy the transformation properties
\begin{equation}
\label{eqn:sl2zsiegel}
 {_c}g_{v \gamma}   = {_c}g_v | \gamma   \qquad \mbox{ for all } v=(\alpha,\beta) \in (\Q/\Z)^2, \quad \gamma \in \SL_2(\Z).
 \end{equation}
(Cf.~\cite[Lemma 1.7(1)]{kato}.)
In particular, ${_c}g_{\alpha,\beta}$ is a unit on the open modular curve attached to the congruence subgroup of $\SL_2(\Z)$ that fixes $(\alpha,\beta)$, and hence belongs to 
$\cO^\times(Y_0(N))$. 
The Siegel units also satisfy the distribution relations:
\begin{equation}
\label{eqn:dist-rel-siegel}
 \prod_{m\alpha'=\alpha}  {_c}g_{\alpha',\beta}(z)  =  {_c}g_{\alpha,\beta}(z/m), 
\qquad \prod_{m\beta'= \beta} {_c} g_{\alpha,\beta'}(z) = {_c}g_{\alpha,\beta}(mz),
\end{equation}
which together  imply that
\begin{equation}
\label{eqn:full-dist-siegel}
 \prod_{m(\alpha',\beta') = (\alpha,\beta)} {_c}g_{\alpha',\beta'}(z)  = {_c}g_{\alpha,\beta}(z).
 \end{equation}
(Cf.~\cite[Lemma 1.7(2)]{kato} or \cite[Prop.~2.2.1 and 2.2.2]{llz}.)

The unit $_c g_{\alpha,\beta} $ is equal to $ g_{\alpha,\beta}^{c^2} \cdot g_{c\alpha,c\beta}^{-1}$,
where the $q$-expansion of 
$g_{\alpha,\beta} \in \cO^\times(Y(N))\otimes \Q$ is given by
\begin{equation}
\label{eqn:siegel-qexp}
 g_{\alpha,\beta}(q) = -q^{w} 
\prod_{n\ge 0}  (1- q^{n+\alpha} e^{2\pi i \beta}) 
\prod_{n>0} (1-q^{n-\alpha} e^{-2\pi i \beta}),  
\end{equation}
where $
w = 1/12 - \alpha/2 + (1/2)\alpha/N$, with $0\le \alpha<1$. 
 (Cf.~\cite[\S 1.9]{kato}.) 

Fix a rational prime $p$, and 
assume that $(\alpha,\beta)$ is of $p$-power order in $(\Q/\Z)^2$. 
 To  lighten notations, it will be assumed below that $p\ne 5$,
  and the choice
 $c=5$ will be fixed. (The constructions  are readily adapted to the case
   $p=5$ by  changing the value of  $c$.)

\subsection{The Siegel distribution}

Let
 $\bX_0:= (\Z_p^2)'$ be the set of vectors $(a,b)\in \Z_p^2$ that are primitive, i.e., satisfy
$\gcd(a,b)=1$, and let 
\begin{equation}
\label{eqn:X-X0}
 \bX := (\Q_p^2-\{0,0\}) = \bigcup_{j=-\infty}^\infty p^j \bX_0.
 \end{equation}
 
 Let $A$ be an  $\SL_2(\Z)$-module, and let 
$\LC(\bX_0,\Z)$ be the space of locally constant $\Z$-valued functions on 
$\bX_0$. An {\em $A$-valued distribution} on $\bX_0$ is a homomorphism from $\LC(\bX_0,\Z)$ to $A$. Because $\bX_0$ is compact, a distribution $\mu$  is determined by its values $\mu(U)$ on the characteristic functions of compact open subsets $U\subset \bX_0$. Let $\bD(\bX_0,A)$ denote the module of $A$-valued distributions. It is endowed with the (right)  $\SL_2(\Z)$-action defined by
\begin{equation}
\label{eqn:action-on-distributions}
 (\mu|\gamma)(U) = \mu(U\gamma^{-1})|\gamma, \qquad \mbox{ for } \gamma\in \SL_2(\Z), \ \ U\subset \bX_0.
 \end{equation}
 A distribution on $\bX$ is said to be $p$-invariant if it is
  invariant under multiplication by $p$,
 i.e., 
 \begin{equation}
 \mu(p^jU) = \mu(U) \mbox{ for all  $j\in \Z$ and all compact open } U \subset \bX.
 \end{equation}
 Denote by $\bD(\bX,A)$   the module of $p$-invariant distributions  on $\bX$.
  Because 
 $\bX_0$ is a fundamental domain for the action of $p$ on $\bX$ (cf.~\eqref{eqn:X-X0}),
 every distribution on $\bX_0$ extends uniquely to a 
 $p$-invariant distribution, yielding  an isomorphism
\begin{equation}
\label{eqn:ident-X0-X}
\bD(\bX_0,A) \stackrel{\simeq}{\lra} \bD(\bX,A).
\end{equation}
 The target space is equipped with a natural action of 
 the larger group $\Gamma$ when $A$ is a $\Gamma$-module, defined by 
 \eqref{eqn:action-on-distributions} with $\SL_2(\Z)$ replaced by $\Gamma$. 
 For all
  $\mu\in \bD(\bX, A)$ and  for all locally constant, compactly supported $\Z$-valued
 functions $f$ on $\bX$, the $\Gamma$-action is determined  by
$$ \int_{\bX} f(x,y) d(\mu|\gamma)(x,y) =  \int_{\bX} f( (x,y)\gamma) d\mu(x,y).$$

As was implicitly observed in the work of Kubert and Lang, the collection of 
Siegel units of $p$-power level  are conveniently packaged into a distribution on 
$\bX_0$, by setting 
$$ \mu_{\rm Siegel}\left((a,b)+p^n(\Z_p^2)\right) := _cg_{\frac{a}{p^n},\frac{b}{p^n}},  \quad
\mbox{ for all } (a,b) \in (\Z^2)'.$$
Since every compact open subset of $\bX_0$ is a union of sets of the form
$(a,b)+p^n(\Z_p^2)$, the above rule determines $\mu_{\rm Siegel}$  on all compact open subsets of $\bX_0$. The fact that it is well-defined follows from the distribution relation
\eqref{eqn:full-dist-siegel} with $m=p$.

View $\mu_{\rm Siegel}$ as an element of $\bD(\bX,\cO_{\cH}^\times)$ via \eqref{eqn:ident-X0-X}. A key feature of  $\mu_{\rm Siegel}$ is its invariance under  $\Gamma=\SL_2(\Z[1/p])$, and even under the full group $\GL_2^+(\Z[1/p])$ of invertible matrices 
 with coefficients in $\Z[1/p]$ and positive determinant. 
\begin{theorem}
\label{thm:kubert-lang}
The distribution $\mu_{\rm Siegel}$ satisfies 
\begin{equation}
\label{eqn:equiv-siegel}
 \mu_{\rm Siegel}(U \gamma) = \mu_{\rm Siegel}(U)|\gamma, 
\end{equation}
for all compact open subsets $U\subset \bX$ and all $\gamma\in \GL_2^+(\Z[1/p])$. 
\begin{proof}
Let $(\alpha,\beta) = (\frac{a}{p^n}, \frac{b}{p^n})$ be an element of order $p^n$ in $(\Q/\Z)^2$.
 Since the sets $U_{\alpha,\beta} = (a,b) + p^n\Z_p^2$ and their translates under multiplication by $p$ form a basis for the topology on $\bX$, it suffices to prove the theorem for the sets of this form. The equivariance \eqref{eqn:equiv-siegel} for $\gamma\in \SL_2(\Z)$ follows directly from \eqref{eqn:sl2zsiegel}. Since $\GL_2^+(\Z[1/p])$ is generated by $\SL_2(\Z)$ and the matrix $T:= \smallmat{p}{0}{0}{1}$, one is reduced to showing the relation 
\[
\mu_{\rm Siegel}(U_{\alpha,\beta} T) = \mu_{\rm Siegel}(U_{\alpha,\beta})|T.
\]
To see this, note that
\begin{eqnarray*}
 U_{\alpha,\beta} T 
 &=&  (pa+ p^{n+1} \Z_p) \times (b+p^n \Z_p)  \\
 &=&  \bigcup_{b'\equiv b  (p^n)} (pa+ p^{n+1} \Z_p) \times (b'+p^{n+1} \Z_p) 
 = \bigcup_{p\beta' 
 = \beta} U_{\alpha, \beta'}.
\end{eqnarray*}
It then follows from \eqref{eqn:dist-rel-siegel} that
\[
\mu_{\rm Siegel}(U_{\alpha,\beta} T) 
= \prod_{p\beta'=\beta} {_cg_{\alpha,\beta'}}(z) 
= \  _cg_{\alpha,\beta}(pz)  
= \mu_{\rm Siegel}(U_{\alpha,\beta})| T,
\]
 as was to be shown.
\end{proof}
\end{theorem}

The invariance of $\mu_{\rm Siegel}$ under translation by the full $p$-arithmetic 
group $\Gamma$, which    is hinted at in  
\cite[Rem.~2.2.3]{llz},   combines 
  the  $\SL_2(\Z)$-invariance properties
\eqref{eqn:sl2zsiegel}  and  norm compatibility relations
   \eqref{eqn:dist-rel-siegel},
\eqref{eqn:full-dist-siegel}
    satisfied by
   the Siegel units of $p$-power
 level into a single unified statement.
 
 The following Lemma evaluates the Siegel distribution at some distinguished open subsets of $\bX$. 
 \begin{lemma}
 \label{lemma:eval-mu-siegel}
 The distribution $\mu_{\rm Siegel}$ satisfies
\begin{equation}
\label{eqn:Siegel_on_opens}
\begin{array}{llll}
\mu_{\rm Siegel}(\bX_0) &=& 1 & \pmod{\pm p^{\Z}},\\
\mu_{\rm Siegel}(p\Z_p\! \times\! \Z_p^{\times}) &=& (\Delta(q^p)/\Delta(q))^2 & \pmod{\pm p^{\Z}}.
\end{array}
\end{equation}
 \begin{proof} 
The first assertion follows from the fact that $\bX_0$ is stabilised by $\SL_2(\Z)$, and therefore that its associated Siegel unit is  a unit  on the open modular curve $Y_0(1)$ of level $1$, which contains no non-constant elements. More precisely, $\mu_{\rm Siegel}(\bX_0)$ belongs to $\cO^\times(Y_0(1)_{\Z[1/p]}) = \pm p^{\Z}$.
(Cf.~\cite[Prop.2.3.2]{llz}  for instance.)
 The second assertion follows from the calculation
 $$ \mu_{\rm Siegel}(p\Z_p\times \Z_p^\times) = \prod_{i=1}^{p-1} \mu_{\rm Siegel}((0,i) + p\Z_p^2) = \prod_{i=1}^{p-1} {_cg}_{0,i/p}  =
 \pm p^{c^2-1} (\Delta(q^p)/\Delta(q))^{(c^2-1)/12},$$
 where the last equality can be read off from the $q$-expansions of the Siegel units given in \eqref{eqn:siegel-qexp}. The result now follows, since 
 $c=5$.
  \end{proof}
 \end{lemma}

\subsection{The Dedekind--Rademacher distributions}
The following general Lemmas concerning $p$-invariant distributions will be useful later.
\begin{lemma}
 \label{lemma:simple-but-important}
 Let $\mu$ be any element of $\bD(\bX,A)$. 
 If $\Lambda$ is any $\Z_p$-lattice in $\Q_p^2$, and $\Lambda'$ is its set of 
 primitive vectors, then 
 $ \mu(\Lambda') = \mu(\bX_0)$. 
\begin{proof}
By compactness, there is an integer $N\ge 0$ for which $p^{-N}\Z_p^2 \subset \Lambda \subset p^N \Z_p^2$, and hence each 
$v\in \Lambda'$ belongs to a  translate $p^j \bX_0$ for a unique $j\in [-N,N]$.  Hence one may write 
$$ \Lambda' = p^{m_1} U_1 \sqcup \cdots  \sqcup p^{m_t}  U_t,$$
for a suitable decomposition 
$$ \bX_0 = U_1 \sqcup \cdots \sqcup U_t$$
of $\bX_0$ as a disjoint union of compact open subsets.
The additivity properties   of $\mu$ combined with its  $p$-invariance    implies  that
$\mu(\Lambda') = \mu(\bX_0)$, as claimed.
\end{proof}
\end{lemma}

\begin{lemma}
\label{lemma:is-exact}
The rule which to $A$ associates  $\bD(\bX,A)$ is 
  an exact (covariant) functor from the category
of $\Gamma$-modules to itself. 
\begin{proof}
The issue is right exactness. If $\varphi: A\lra B$ is a surjective module  homomorphism and 
$\mu\in \bD(\bX,B)$ is a $B$-valued, $p$-invariant distribution on $\bX$,   one can construct
a distribution $\tilde \mu\in \bD(\bX,A)$ that maps to it by
choosing, for  each successsive  $n\ge 1$ and for each primitive vector $v = (\Z/p^n\Z)'$, the value 
$\tilde\mu(v + p^n\Z_p^2) \in A $  satisfying $\varphi(\tilde\mu(v+p^n\Z_p^2)) = \mu(v+p^n\Z_p^2)$, 
taking care at each stage that the additivity relations required of distributions be satisfied. 
One obtains in this way an element of $\bD(\bX_0,B)$, giving rise to the desired lift 
in $\bD(\bX,B)$ via
\eqref{eqn:ident-X0-X}.
\end{proof}
\end{lemma}

Thanks  to Lemma \ref{lemma:is-exact}, the exponential sequence
$$
\xymatrix{ 0 \ar[r]  & \Z  \ar[r] & \cO_{\cH} \ar[r]^{ e^{2\pi i z}} &  {\cO_{\cH}^\times} \ar[r] & 1 }$$
 induces a short exact sequence
$$ 1 \rightarrow \bD(\bX,\Z) \lra \bD(\bX, \cO_{\cH}) \lra \bD(\bX, \cO_{\cH}^\times) \rightarrow 1$$
 of
$\Gamma$-modules.
Let 
$$ \delta: \rH^0(\Gamma, \bD(\bX, \cO_{\cH}^\times)) \lra \rH^1(\Gamma, \bD(\bX,\Z))$$
be the connecting homomorphism arising from the 
resulting long exact sequence in $\Gamma$-cohomology.
The image 
$$\mu_{\rm DR} := \delta(\mu_{\rm Siegel})  \in \rH^1(\Gamma,\bD(\bX, \Z))$$ 
is a one-cococyle on $\Gamma$, i.e., it satisfies the relation
$$ \mu_{\rm DR}(\gamma_1\gamma_2) = \mu_{\rm DR}(\gamma_1)   + \mu_{\rm DR}(\gamma_2)| \gamma_1^{-1}.$$
It is obtained by lifting $\mu_{\rm Siegel}$ to an $\cO_{\cH}$-valued distribution 
$$\tilde \mu_{\rm Siegel}  := \frac{1}{2\pi i} \log(\mu_{\rm Siegel}) \in \bD(\bX,\cO_{\cH}),$$
and setting
\begin{equation}
\label{eqn:tilde-siegel}
 \mu_{\rm DR}(\gamma) := \tilde\mu_{\rm Siegel}|{\gamma^{-1}}  - \tilde \mu_{\rm Siegel}.
 \end{equation}

Recall the Dedekind--Rademacher homomorphism $\varphi_{\rm DR}: \Gamma_0(p) \lra \Z$  evoked in the introduction, which encodes the periods of the  Eisenstein series $E_2^{(p)} = d\log(\Delta(pz)/\Delta(z))$ of weight two. 
\begin{lemma} 
\label{lemma:eval-muDR}
The one-cocycle $\mu_{\rm DR}$ satisfies
\begin{equation}
\begin{array}{llll}
\mu_{\rm DR}(\gamma)(\bX_0) &=& 0 & \qquad \mbox{for all } \ \gamma\in \Gamma, \\
\mu_{\rm DR}(\gamma) (p\Z_p\! \times\! \Z_p^\times) &=&\varphi_{\rm DR}(\gamma) & \qquad \mbox{for all } \ \gamma\in \Gamma_0(p). 
\end{array}
\end{equation}
 \begin{proof} 
Observe that, for all $\gamma\in \Gamma$,
 $$ \mu_{\rm DR}(\gamma)(\bX_0) =  \tilde\mu_{\rm Siegel}|{\gamma^{-1}}(\bX_0)  - \tilde \mu_{\rm Siegel}(\bX_0) = \tilde\mu_{\rm Siegel}(\bX_0\gamma)|\gamma^{-1} - \tilde\mu_{\rm Siegel}(\bX_0).$$
 Lemma   \eqref{lemma:simple-but-important}  implies that $ \tilde\mu_{\rm Siegel}(\bX_0\gamma) =\tilde\mu_{\rm Siegel}(\bX_0)$, and 
 Lemma \ref{lemma:eval-mu-siegel} shows that  this common value is a constant function on $\cH$.  The first assertion follows.
 As for the second,
 equation \eqref{eqn:tilde-siegel} implies that
 $$ \mu_{\rm DR}(\gamma) (p\Z_p\times \Z_p^\times) = (\tilde\mu_{\rm Siegel}|{\gamma^{-1}}  - \tilde \mu_{\rm Siegel})(p\Z_p\times \Z_p^\times).$$
 By Lemma \ref{lemma:eval-mu-siegel}, 
 $$  \tilde \mu_{\rm Siegel}(p\Z_p\times \Z_p^\times) = \frac{2}{2\pi i} \log\left(\Delta(pz)/\Delta(z)\right) \pmod{\C}.$$
 Since $\gamma\in \Gamma_0(p)$ preserves the region
 $p\Z_p\times\Z_p^\times$, it follows that
 $$  (\tilde\mu_{\rm Siegel}|{\gamma^{-1}}  - \tilde \mu_{\rm Siegel})(p\Z_p\times \Z_p^\times) = \frac{2}{2\pi i} \int_{z_0}^{\gamma^{-1} z_0} d\log(\Delta(pz)/\Delta(z)) 
 = \varphi_{\rm DR}(\gamma),$$
 as was to be shown.
  \end{proof}
  \end{lemma}

\subsection{The multiplicative Poisson transform}

Because a  distribution  $\mu\in \bD(\bX,\Z)$ is $\Z$-valued, and hence $p$-adically bounded, it also gives rise to a
{\em measure}:  one can extend 
$ \mu$ to arbitrary continuous, compactly supported functions on $\bX$.
There is even a {\em multiplicative refinement} of the integral against $\mu$, defined by
$$ \mint_{\bX} f(x,y) d\mu(x,y) := \lim_{\{U_\alpha\}} \prod_\alpha f(x_\alpha,y_\alpha)^{\mu(U_\alpha)},$$
where the limit is taken over finer and finer open covers $\{U_\alpha\}$ of the support of $f$, and
$(x_\alpha,y_\alpha)$ is a sample point in $U_\alpha$.  Here
 $f:\bX\lra \C_p^\times$ is a continuous, compactly supported function on $\bX$ (which means that it takes the value $1$
outside a compact subset of $\bX$).

Let $\bD_0(\bX_0,\Z)$ be the $\Z$-module of distributions on $\bX_0$ satisfying
$$ \mu( \bX_0)= 0.$$
The multiplicative Poisson transform of   $\mu\in \bD_0(\bX_0,\Z)$ is the rigid analytic function $J(\mu)$  on $\cH_p$ 
defined by setting
$$ J(\mu)(\tau) = \mint_{\bX_0} (x\tau + y) d\mu(x,y).$$
This assignment gives rise to an $\SL_2(\Z)$-equivariant map
$$ J: \bD_0(\bX_0,\Z) \lra \cA^\times,$$ i.e., 
$$ J({\mu|\gamma})(\tau) = J(\mu) | \gamma(\tau)= J(\mu)(\gamma\tau), 
\qquad \mbox{ for all } \gamma \in \SL_2(\Z). $$
Identifying $\bD_0(\bX_0,\Z)$ with the module $\bD_0(\bX,\Z)$  
of distributions on $\bX$  satisfying
$$ \mu( \bX_0)= 0, \qquad \mu( pU) = \mu(U), $$
the same rule $J$ (where one continues to integrate over the compact subset $\bX_0\subset \bX$)
determines  a $\Gamma$-equivariant map
\begin{equation}
\label{eqn:full-inv}
 J: \bD_0(\bX,\Z) \lra \cA^\times/p^{\Z}.
 \end{equation}
 The reason for this somewhat weaker invariance property is that while $\SL_2(\Z)$ preserves the region
$\bX_0$ of integration defining $J(\mu)$, the full $p$-arithmetic  group $\Gamma$ does not.
Nonetheless, if $\gamma\in \Gamma$, one still can write  (following the reasoning in the proof of
  \ref{lemma:simple-but-important})
$$\bX_0\gamma = p^{m_1} U_1 \sqcup \cdots  \sqcup p^{m_t} U_t, \qquad \mbox{ with } 
\quad \bX_0 = U_1 \sqcup \cdots \sqcup U_t,$$
 and  the integrand
$(x-\tau y)$ arising in the definition of $J$  obeys a simple transformation property under multiplication by $p$.
It follows that $J(\mu|\gamma) = J(\mu)|\gamma \pmod{p^{\Z}}$, for all $\gamma\in \Gamma$.

Let 
 $$J_{\rm DR} := J(\mu_{\rm DR}) \in \rH^1(\Gamma,\cA^\times/p^{\Z})$$
  be the image of
  the measure-valued cocycle $\mu_{\rm DR}$ under the multiplicative
   Poisson transform of \eqref{eqn:full-inv}.
It is represented by the  one-cochain $J_{\rm DR}: \Gamma\lra \cA^{\times}$ (denoted by the same symbol, by an abuse of notation)
defined by
   $$ J_{\rm DR}(\gamma)(\tau) = J(\mu_{\rm DR}(\gamma))(\tau),$$
  which satisfies the  cocycle relation modulo $p^{\Z}$,
  $$ J_{\rm DR}(\gamma_1 \gamma_2) = J_{\rm DR}(\gamma_1) \times J_{\rm DR}(\gamma_2) | \gamma_1^{-1} \pmod{p^{\Z}}.$$
 Its restriction to $\SL_2(\Z)$ also satisfies the  {\em full cocycle 
 relation}, with no $p^{\Z}$-ambiguity, because of the $\SL_2(\Z)$-equivariance of $J$.

\medskip
In order to prove Theorem A of the introduction, it now suffices to
calculate the image of $J_{\rm DR}$ under the sequence of maps
$$ \eta:  \rH^1(\Gamma, \cA^\times/p^{\Z}) \lra \rH^2(\Gamma, p^{\Z}) = \rH^1(\Gamma_0(p), p^{\Z}).$$
\begin{theorem} 
The image of $J_{\rm DR}$  under $\eta$ is
$$ \eta(J_{\rm DR}) = p^{\varphi_{\rm DR}}.   $$
\begin{proof}
The action of $\Gamma$ on the Bruhat-Tits tree of $\PGL_2(\Q_p)$ leads to an expression for $\Gamma$ as an amalgamated product of the groups
$$ \SL_2(\Z), \qquad \SL_2(\Z)' = \mat{p}{0}{0}{1}^{-1} \SL_2(\Z) \mat{p}{0}{0}{1},$$
whose intersection is $\Gamma_0(p)$. 
The fact that $\rH^1(\SL_2(\Z), p^{\Z}) = 0$
and that $\rH^2(\SL_2(\Z),p^{\Z})$ is of order $12$  ensures the existence
of unique lifts to $\cA^\times$ of the restrictions of $J_{\rm DR}^{12}$ to $\SL_2(\Z)$ and
$\SL_2(\Z)'$:
$$\JJ_{\rm DR}  \in \rH^1(\SL_2(\Z), \cA^{\times}),  \qquad \JJ_{\rm DR}' \in \rH^1(\SL_2(\Z)',\cA^{\times}).$$
One  then has, for all $\gamma\in \Gamma_0(p)$, 
\begin{equation}
\label{eqn:formule-eta}
 \eta(J_{\rm DR}^{12})(\gamma) = \JJ_{\rm DR}(\gamma) \div \JJ_{\rm DR}'(\gamma).
 \end{equation}
Concretely,   $\JJ_{\rm DR}$ 
and $\JJ_{\rm DR}'$ may be expressed
as multiplicative Poisson transforms of $\mu_{\rm DR}$,
by setting
$$ \JJ_{\rm DR}(\gamma)(\tau):= \mint_{\bX_0} (x\tau + y)^{12} d\mu_{\rm DR}(\gamma)(x,y), 
\qquad 
\JJ_{\rm DR}'(\gamma)(\tau) :=  \mint_{\bX_0'} (x\tau + y)^{12} d\mu_{\rm DR}(\gamma)(x,y),
$$
where $\bX_0' := (p\Z_p \times \Z_p)'$ is the translate of $\bX_0$ under the matrix
$\mat{p}{0}{0}{1}$, a region whose stabiliser in $\Gamma$ is the group $\SL_2(\Z)'$. 
Observe that
\begin{equation}
\label{eqn:overlaps}
 \bX_0 \cap \bX_0' = p\Z_p\times \Z_p^{\times}, \qquad
 \bX_0 - \bX_0' = \Z_p^{\times} \times \Z_p, \qquad
 \bX_0'-\bX_0 = p(\Z_p^{\times} \times \Z_p).
 \end{equation}
 Hence, for all $\gamma\in \Gamma_0(p)$,
\begin{eqnarray*}
 \JJ_{\rm DR}(\gamma) \div \JJ_{\rm DR}'(\gamma)  &=& 
 \mint_{\bX_0} (x\tau + y)^{12} d\mu_{\rm DR}(\gamma)(x,y) \div \mint_{\bX_0'} (x\tau + y)^{12} d\mu_{\rm DR}(\gamma)(x,y) \\
 &=& \mint_{\Z_p^{\times}\times \Z_p} (x\tau + y)^{12} d\mu_{\rm DR}(\gamma)(x,y) \div \mint_{p(\Z_p^{\times}\times\Z_p)} (x\tau + y)^{12} d\mu_{\rm DR}(\gamma)(x,y) \\
 &=& \mint_{ \Z_p^{\times}\times \Z_p} p^{-12} d\mu_{\rm DR}(\gamma)(x,y),
 \end{eqnarray*}
 where the penultimate equality follows from  \eqref{eqn:overlaps} and the 
 last from   the invariance of $\mu_{\rm DR}(\gamma)$ under
 multiplication by $p$.
 Because $(\Z_p^\times \times \Z_p)$ is the complement of $(p\Z_p\times \Z_p^{\times})$ in $\bX_0$, 
 and $\mu_{\rm DR}(\gamma)(\bX_0)=0$, this implies that
$$  \JJ_{\rm DR}(\gamma) \div \JJ_{\rm DR}'(\gamma)  = \mint_{p\Z_p\times\Z_p^{\times}} p^{12} d\mu_{\rm DR}(\gamma)(x,y) = 
p^{12\mu_{\rm DR}(\gamma)(p\Z_p\times\Z_p^{\times})} = p^{12\varphi_{\rm DR}(\gamma)},$$
where the last equality follows from Lemma \ref{lemma:eval-muDR}. Combining this with
\eqref{eqn:formule-eta}
shows that $\eta(J_{\rm DR})$ and $p^{\varphi_{\rm DR}}$ agree, since the group they belong to is torsion-free.
This completes the proof of Theorem A. 
\end{proof}
\end{theorem}

\section{The winding cocycle}
\label{sec:winding}

The goal of this section is to recall the definition and key properties of the winding cocycle introduced in \cite[\S~2.3]{DPV1} and to  reduce Theorem B of the introduction to Theorem C.

\subsection{The residue map}
The group $
\rH^1(\Gamma, \cA^\times/\C_p^\times)
$
of rigid analytic theta cocycles
is finitely generated and closely related to the space of modular forms of weight two on the Hecke congruence group $\Gamma_0(p)$. More precisely, it is a module over the Hecke algebra $\TT_0(p)$ of Hecke operators acting faithfully on the weight two modular forms on $\Gamma_0(p)$. To see this, let 
\[
U := \{ z\in \PP_1(\C_p) \mbox{ with } 1 < |z| < p \} \subset \cH_p
\]
be the standard annulus whose stabiliser in $\Gamma$ is $\Gamma_0(p)$. The logarithmic annular residue map 
   \begin{equation}
   \label{eqn:def-residue-V}
    \partial_{U}: \cA^\times/\C_p^\times \lra \Z_p, \qquad  \partial_{U}(f) := {\rm Res}_{U}(\dlog f)
    \end{equation}
is equivariant for the action of $\Gamma_0(p)$, and hence composing it with the restriction to $\Gamma_0(p)$ yields a map on cohomology
     \begin{equation}
   \label{eqn:def-residue-V-coh}
   \partial_{U}: \rH^1(\Gamma, \cA^{\times}\! /\C_p^{\times}) \lra 
   \rH^1(\Gamma_0(p),\Z_p),
   \end{equation}
which is denoted by the same symbol by abuse of notation. This map is compatible with the action of the Hecke operators, and with the involution $W_\infty$ determined by the matrix 
 $\smallmat{1}{0}{0}{-1}$, which lies in the normaliser of both $\Gamma$ and $\Gamma_0(p)$.
 Let $\rH^1(\Gamma, \cA^{\times}\! /\C_p^{\times})^\pm$ denote the plus and minus eigenspaces for this involution in the space of rigid analytic theta cocycles, and  denote by $\rH^1(\Gamma_0(p),\Z_p)^\pm$ the corresponding eigenspaces in the cohomology of $\Gamma_0(p)$.
 
\par While the map in \eqref{eqn:def-residue-V} has an infinite rank  kernel, it is notable that the induced map on rigid analytic theta cocycles is essentially an isomorphism:
  \begin{lemma} 
   \label{lemma:residue-isomorphism}
      Up to torsion kernels and cokernels,
      the map $\partial_{U}$ of  
        \eqref{eqn:def-residue-V-coh}
      is surjective, and its kernel is generated 
      by the ``trivial" theta-cocycle 
      \[
      J_{\rm triv} \in \rH^1(\Gamma, \cA^\times), \qquad 
      J_{\rm triv}       \left(\smallmat{a}{b}{c}{d}\right)(z) = cz+d .
      \]
      In particular, the induced map
    \begin{equation}
   \label{eqn:def-residue-V-coh-minus}
   \partial_{U}^- \ : \   \Q \otimes \rH^1(\Gamma, \cA^{\times}/\C_p^{\times})^- \ \lra  \ 
   \rH^1(\Gamma_0(p),\Q)^-
   \end{equation}
 is an isomorphism.
 \begin{proof}
The first assertion is  a reformulation of \cite[Theorem 3.1]{DPV1}. The last follows from the fact that $J_{\rm triv}$ is fixed by  $W_\infty$, as can be checked directly from the definition of $J_{\rm triv}$.
\end{proof}
      \end{lemma}

  \subsection{The winding cocycle}
  \label{sec:winding-cocycle}
  In \cite[\S 2.3]{DPV1}, the so-called 
  {\em winding cocycle}
    $$ J_w \in \rH^1(\Gamma, \cA^\times/\C_p^\times)^-$$
  is introduced. Unlike the Dedekind--Rademacher cocycle, it is not an eigenclass for the Hecke operators, although it belongs to the $-1$ eigenspace for the involution $W_\infty$. The greater complexity of $J_w$ on the spectral side is offset by a gain in simplicity on the geometric side, evidenced by  the fact that the rigid analytic functions $J_w(\gamma)$ admit explicit infinite product expansions. 
  
\par Let $Y_0(p) = \Gamma_0(p) \backslash \cH$ be the open modular curve, and let $X_0(p)$ be its standard compactification, obtained by adding the two cusps $0$ and $\infty$. The intersection pairing on homology (Poincar\'e duality) 
  defines  isomorphisms
  \begin{equation}
  \label{eqn:pd-iso}
\rH_1(X_0(p) ; \{0,\infty \}, \Q)^\pm  = \rH^1(Y_0(p), \Q)^\mp
 = \rH^1(\Gamma_0(p), \Q)^\mp.
\end{equation}
Mazur's {\em winding element}    
\[
\varphi_w \in  \rH^1(\Gamma_0(p),\Z)^-  
\]
is defined to be the class of the path from $0$ to $\infty$ in the  homology of the modular curve $X_0(p)$ relative to the cusps, viewed as an element of $\rH^1(\Gamma_0(p),\Z)$ via \eqref{eqn:pd-iso}. By \cite[Prop.~3.3]{DPV1} and its proof, the winding cocycle is characterised by the identity  
\begin{equation}
\label{eqn:residue-Jw}
 \partial_U^-(J_w) = 2 \varphi_w.
 \end{equation}

  \subsection{Theorem C implies Theorem B}
  
Theorem B of the introduction is  reduced to Theorem C   by  writing
 the modular form $G_\tau$  of this theorem
 as a linear  combination of eigenforms.
 
To this end,   observe that
 $\rH^1(\Gamma_0(p), \bar\Q)^-$ is generated as a $\bar\Q$-vector space by 
  the Dedekind--Rademacher morphism 
  $\varphi_{\rm DR}$ of
  \eqref{eqn:varphi-DR}
  encoding the periods of the weight two Eisenstein series $E_2^{(p)}$ defined in the introduction, 
  and the homomorphisms $\varphi_f^-$ attached to the  minus modular symbol  for $f$, where $f$ runs through a basis of cuspidal Hecke eigenforms in $S_2(\Gamma_0(p))$. 
  A direct calculation of  integration pairings in \cite[Lemma 3.4]{DPV1} then yields the  spectral decomposition of the winding element,
\begin{equation}
\label{eqn:hom-spectral}
\varphi_w \ = \
 \frac{1}{p-1}\cdot\varphi_{\rm DR} \ + \ \sum_f \lambda_f \cdot \varphi_f^-,
\end{equation}
where the coefficient $\lambda_f\in \bar\Q$ is a suitable non-zero
multiple of $L(f,1)$  whose exact nature is not germane to the proof of Theorem B. (But see \cite[\S~3]{DPV1} for more details.)

\par Now consider the rigid analytic theta-cocycle
\[
J_f^- \ \in \   \bar \Q \otimes \rH^1(\Gamma,\cA^\times/\C_p^\times)^-
\]
characterised by $\partial_U^-(J_f^-) = \varphi_f^-$. By Lemma \ref{lemma:residue-isomorphism} and \eqref{eqn:residue-Jw}, 
\begin{equation}
\label{eqn:ratc-spectral}
J_w \ = \ \frac{2}{p-1} \cdot J_{\rm DR}   \ \ + \ \  \sum_f  2\lambda_f \cdot J_f^-  \qquad \mbox{in } \  \bar\Q \otimes \rH^1(\Gamma,\cA^\times/\C_p^\times)^-,
\end{equation}
 where additive notation has been adopted to describe the operations in this group in spite of its multiplicative nature.
For each $n\ge 1$, applying the Hecke operator $T_n$ to this identity
then gives 
\begin{eqnarray}
\label{eqn:wind-spectral}
 T_nJ_w &=&  \frac{2}{p-1} \cdot T_n J_{\rm DR}   \ \ + \ \ 
 \sum_f  2\lambda_f \cdot  T_nJ_f^-  \\
 \nonumber
  &=&   \frac{2}{p-1}  \cdot  J_{\rm DR} \cdot \sigma_1^{(p)}(n)     \ \ + \ \  \sum_f  2\lambda_f  \cdot  J_f^- \cdot a_n(f)  
\end{eqnarray}
in $\bar\Q \otimes \rH^1(\Gamma,\cA^\times/\C_p^\times)^-$. After evaluating at the RM point $\tau$ and taking $p$-adic logarithms, it follows that
\begin{equation}
\label{eqn:wind-spectral-bis}
 \log_p(T_nJ_w[\tau]) \ \ \ =  \ \ \ \frac{2 \log_p(J_{\rm DR}[\tau])}{p-1} \cdot\sigma_1^{(p)}(n)  \ \   + \ \ \sum_f  2 \lambda_f   \log_p(J_f^-[\tau]) \cdot a_n(f). 
\end{equation}
 Substituting this identity into 
 Theorem C of the introduction yields the spectral expansion 
\begin{equation}
\label{eqn:spectral-decomp-Gtau}
 G_\tau(q)  \ \ =  \ \  \frac{ \log J_{\rm DR}[\tau]}{ 12(p-1)}  \cdot E_2^{(p)}(q) \ +  \ 
  \sum_{f} \beta_f \cdot f(q),
  \end{equation}
  where $E_2^{(p)}$ is the Eisenstein series of 
  \eqref{eqn:E2p}, and 
  \begin{equation}
  \label{eqn:betaf}
  \beta_f = 2 \lambda_f   \log_p(J_f^-[\tau]).
  \end{equation}
 
    Comparing  the zero-th Fourier coefficient
    of $G_\tau$   in     \eqref{eqn:spectral-decomp-Gtau}
    with the one in Theorem C
shows that
   $$ \frac{ \log J_{\rm DR}[\tau]}{12} =   \log(u_\tau),$$
thereby reducing Theorem B of the introduction to Theorem C.

The remainder of the paper is devoted to the construction of the modular generating series required for the proof of Theorem C.

\medskip\noindent
{\bf Remark}. 
The coefficients $\beta_f$ in \eqref{eqn:betaf} are immaterial to the proof of Theorem B but
are of independent interest, insofar as they involve  the RM values of the elliptic rigid analytic theta-cocycles $J_f^-$: these values are the formal group logarithms of certain {\em Stark--Heegner points} in the modular Jacobian $J_0(p)$. Although poorly understood theoretically, these Stark--Heegner points are conjectured to be defined over the narrow ring class field $H_\tau$. 
The approach to the algebraicity of $J_{\rm DR}[\tau]$ based on deformations of Galois representation does not seem to shed any immediate light on the algebraicity  of these more mysterious invariants.

\section{Deformations of Hilbert Eisenstein series}
\label{sec:deformation}

This section studies the  derivatives of certain $p$-adic analytic families
of Hilbert modular forms for $F$ parametrised by the weight and specialising to a certain Hilbert Eisenstein series of parallel weight one. 

\par This Eisenstein series has several notable features. Firstly, it is cuspidal when viewed as a $p$-adic modular form, and admits cuspidal $p$-adic deformations. Secondly, it vanishes upon diagonal restriction. This implies that the derivatives of both cuspidal and Eisenstein families specialising to $f$, in spite of not displaying any simple modularity properties themselves, yield $p$-adic modular forms after taking diagonal restriction. A suitable linear combination of these derivatives is  considered in \S~\ref{sec:generating-series}, and the  Fourier coefficients of its ordinary projection are related to the RM values of the winding cocycle.

\par While $p$-adic Eisenstein families only occur in
parallel weight, cuspidal families vary over a larger weight space. The main result of this section is Theorem \ref{thm:alice}, which describes the Fourier coefficients of the derivatives of a cuspidal family in the ``anti-parallel'' direction of the weight space. Much like in the archimedean settings, the Fourier expansions of $p$-adic Eisenstein families are entirely explicit; however, no general expression is available for cuspidal families.  
Our approach to studying cuspidal  deformations of a Hilbert Eisenstein series rests   on the analysis of  the associated Galois  deformation problems. Roughly speaking,  first order deformations of the Artin representation attached to  a Hilbert Eisenstein series of parallel weight one are described in terms of the Galois cohomology of the adjoint representation, which cuts out a finite abelian extension $H$ of $F$. A class in the Galois cohomology of the adjoint  cuts out an abelian $p$-adic Lie  extension  of   $H$, and  the Frobenius traces on the associated Galois deformation involve  $p$-adic logarithms of  global $p$-units in $H$,  via the reciprocity law of global class field theory for $H$. This  translates into the appearance of the logarithms of Gross--Stark units in the Fourier coefficients of first order deformations of  Hilbert Eisenstein series, and accounts   for the presence of the same quantities
 in the constant term of the generating series $G_\tau$ of Theorem C. 

The Galois deformation arguments are clarified  and not substantially lengthened  by working in the setting where $F$ is an arbitrary totally real field of degree $d$ in which $p$ is inert. This will be assumed until \S~\ref{subsec:q-exp}, when the main results will be specialised to the case where $F$ is real quadratic.

\subsection{Hilbert modular forms and Hecke algebras.}
\label{subsec:setup}

Fix an embedding $\bar \Q \hookrightarrow \bar \Q_p$. 
Let $F$ be a totally real field in which $p$ is inert, and
denote by $\mathfrak d$  the different 
of  its ring of integers $\cO_{F}$.
   Write $\alpha_1, \ldots, \alpha_d$ for the distinct embeddings of $F$ into $\bar \Q_p$, so that $\alpha_1$ is the embedding given by the restriction to $F$ of the chosen embedding $\bar \Q \hookrightarrow \bar \Q_p$. 
Via the choice of an isomorphism $\C \simeq \C_p$, one obtains a corresponding indexing of embeddings $\bar \Q \hookrightarrow \C.$
For $x\in F$,  let $(x_1, \dots, x_d)$ denote the image of $x$ under the  embeddings $\alpha_1,\ldots,\alpha_d$, viewed
as a $d$-tuple of either complex or $p$-adic numbers, depending on the context.\\
\emph{It is assumed throughout that the Leopoldt Conjecture holds for $F$}. (When $F$ is  quadratic, this assumption is known to be satisfied.)  

 Fix a totally odd character $\psi$ of the narrow class group of $F$.   Let $E$ be a finite extension of $\Q_p$ containing the images  of $F$ under all embeddings $\alpha_1, \ldots, \alpha_d$ and the values of the character $\psi$.

\par We now recall some definitions and conventions
 related to Hilbert modular forms and their associated Hecke algebras, following the treatment that is given in \cite{Shimura}, \cite[\S~2]{Hid88} and \cite[\S~3]{Hid91}. 

Let ${\bf k}=( k_j) \in \Z^d_{\geq 2}$ be a $d$-tuple of integers. Denote  ${\bf t}=( t_j)$ the vector with $ t_j=1$ for every $1 \leq j \leq d$. Choose a vector ${\bf v} \in \Z^d$ of non-negative integers such that ${\bf k} + 2{\bf v} = m{\bf t}$ for some $m \in \Z$, and define $\bf {w} = \bf {k} +\bf {v } -\bf {t}$. 
The space of  {\em Hilbert modular forms} of weight $({\bf k},{\bf w})$ and full level $\fN$, defined as in  \cite[\S~2]{Hid88}, is a finite-dimensional complex vector space.  Let 
$\mathscr{H}_{{\bf k},{\bf w}}(\fN)$ be the algebra of Hecke operators acting faithfully on the subspace of cuspforms. It is free of finite rank as a $\Z$-module.

Fix $({\bf k,\bf w})$ as above. The $p$-adic Hecke algebra is defined to be
\begin{equation}
 \cT \ := \ \varprojlim_{\alpha} \,\mathscr{H}_{\bf{k},\bf{w}}(p^{\alpha})\otimes \cO_{E}, 
\end{equation}
where the inverse limit is taken with respect to restriction of increasing full level structure at $p$. 
 It contains in particular diamond operators $\langle \fl \rangle$ for every integral ideal $\fl$ coprime to $p$, as well as Hecke operators
\begin{equation}
 \mathbf{T}(y) \ = \ \varprojlim_{\alpha} T(y)y_p^{-{\bf v}}
\end{equation}
for any id\`ele $y \in \widehat{\cO_F} \cap \A_F^{\times}$, whose component at $p$ is denoted $y_p$. When $y_p$ is a unit, the operator $T(y)$ depends only on the integral ideal $\fm$ defined by $y$, and we write $T_{\fm}$ and $\mathbf{T}_{\fm}$ for $T(y)$ and $\mathbf{T}(y)$. The compact ring $\cT$ has a unique decomposition $\cT = \cT^{\nord} \oplus \cT^{\scriptstyle \mathrm{ss}}$ such that $\mathbf{T}(p)$ is a unit in $\cT^{\nord}$ and is topologically nilpotent in $\cT^{\scriptstyle \mathrm{ss}}$. The ring $\cT^{\nord}$ is called the  \textit{nearly ordinary cuspidal Hecke algebra}.  It is independent of the choice of $({\bf k}, {\bf w})$. 

 \par Write $\bfU=(\cO_{F}\otimes \Z_p)^\times$ and let $\bfZ$ be the Galois group of the maximal abelian extension of $F$ unramified outside $p$ and $\infty$. The Iwasawa algebra 
\[
 \Lambda  :=\mathcal O_E \lsem \bfU \times \bfZ \rsem     
 \]
is abstractly isomorphic to a ring of power series in several variables with coefficients in a finite group ring over $\cO_E$. Denote 
\[
\kappa^{\univ} \colon \bfU \times \bfZ  \lra   \Lambda^\times
\] 
its universal character.
Since $p$ is inert in $F$, the group
$\bfU$ is identified with the units of  $F_p$.  
Denote $\bfUtor$ and $\bfZtor$ the torsion free parts of $\bfU$ and $\bfZ$ respectively and let $\Lambda^\circ=\cO_E \lsem \bfUtor \times \bfZtor   \rsem$.
Let $\chi_p$ and $\omega_p$ be the cyclotomic and the Teichm\"uller characters of $G_F$ respectively. They factor through the quotient $\bfZ$. Define $\mathsf{q} = 4$ when $p=2$, and $\mathsf{q}=p$ otherwise. Then the homomorphism 
\[
\chi_p \cdot \omega_p^{-1} \colon \bfZ \to 1+\mathsf{q} \Z_p
\]
induces an isomorphism when restricted to $\bfZtor$ if the Leopoldt conjecture holds for $F$. 
For any weight $({\bf k},{\bf w})$ as above, consider the character $\kappa_{{\bf k}, {\bf w}} \colon  \bfG= \bfU \times \bfZ \lra \bar \Q_p^{\times}$ defined by 
\begin{equation}
(a,z) \  \longmapsto \  a^{\bf v} \cdot \chi_p^{m-1}(z).
\end{equation}
With a slight abuse of notation,  the corresponding ring homomorphism will also be denoted by $\kappa_{ {\bf k},{\bf w}}$ and referred to as the weight $(\bfk ,\bfw)$-specialisation. If ${\bfv}=0$, so that $\bfk=m \cdot  \bft$, the pair $(\bfk, \bfw)$ will be called parallel weight $m$. Parallel weight specialisations are parametrised by the Iwasawa algebra $\cO_E \lsem \bfZ \rsem$, which shall be  regarded as a quotient of $\Lambda$.

The nearly ordinary Hecke algebra 
$\cT^{\nord}$ is a $\Lambda$-algebra via the action of the diamond operators. The main theorem of \cite{Hid89a} asserts that  $\cT^{\nord}$ is  finitely generated and torsion-free as a $\Lambda^\circ$-module. In addition,  the quotient of $\cT^\nord$ by the ideal generated by the kernel of $\kappa_{\bfk,\bfw}$, for $(\bfk, \bfw)$ as above,  is isomorphic to the ordinary part of the classical Hecke algebra of weight $(\bfk, \bfw)$ and Iwahori level at $p$.  This result is often referred to as 
{\em Hida's Control Theorem}.

\medskip\noindent
{\em
$p$-adic families of Hilbert modular forms and  weight one Eisenstein series.}

\noindent
The space of Hilbert modular forms of weight $(\bfk, \bfw)$  and any level is automatically cuspidal unless $(\bfk, \bfw)$ is parallel \cite[($1.8_a$)]{Shimura}. However, for parallel weights, non-trivial Eisenstein forms exist and can be interpolated in explicit $p$-adic families parametrised by $\cO_E\lsem \bfZ \rsem$. The study of congruences between cuspidal and Eisenstein families of Hilbert modular forms is at the heart of Wiles' proof of Iwasawa Main Conjecture over totally real fields \cite{wiles}. In a similar spirit, we consider certain cuspidal and Eisenstein families sharing the same specialisation at parallel weight one, i.e. $\bfk = \bft$ and $\bfw = 0$.

 For any pair $(\phi, \eta)$ of unramified characters of $F$ with $\phi \eta$ totally odd, there exists a family $\cE(\phi, \eta)$ with Fourier expansions as described in \cite[\S~3]{DDP11}. The cases where 
 $(\phi,\eta)=(1,\psi)$ or $(\psi,1)$ are of particular relevance in the calculations leading to the proof of Theorem C.

Let $E_1(1, \psi)$ be the classical Eisenstein series of weight 1 and trivial level with Fourier expansion:
\[
E_1(1,\psi)(z_1,\dots,z_d)=L(F, \psi, 0)+\psi^{-1}(\fd)L(F,\psi^{-1},0)+
2^d\sum_{\nu \in {\fd}_+^{-1}}\sigma_{0,  \psi}(\nu \fd)e^{2\pi i (\nu_1 z_1+ \dots+\nu_d z_d)},
\]
where $z=(z_1, \dots z_d)\in \mathcal H^d$ and  
\[
\sigma_{k, \psi}(\alpha)=\sum_{I \lhd \cO_F, I \mid (\alpha)} \psi(I) \mathrm{Nm}(I)^{k} \qquad \mbox{and} \qquad 
L(F, \psi, s)=\sum_{I \lhd \cO_F} \psi(I)\mathrm{Nm} (I)^{-s},
\]
the latter converging for $\mathrm{Re}(s)$ large enough, analytically continued to all $s\in \C$.

\par In the case where $p$ is inert in $F$, one has $\psi(p)= 1$, and the $p$-adic $L$-functions $L_p(\psi,s)$ and $L_p(\psi^{-1},s)$ 
have   exceptional zeros at $s=0$. The Eisenstein series  $E_1(1, \psi)$ then admits a unique $p$-stabilisation 
\begin{equation} \label{p-stabilisation}
f(z):= E_1^{(p)}(1,\psi)(z) := E_{1}(1, \psi)(z)- E_1(1, \psi)(pz)
\end{equation}
which is the weight one specialisation of the Eisenstein families $\cE(1, \psi)$ and $\cE(\psi, 1)$. The derivatives of the Fourier coefficients of $\cE(1, \psi)$ and $\cE(\psi, 1)$ at weight 1 will be exploited in \S \ref{sec:generating-series}. Let 
$$\kappa_{1+\eps} \colon \cO_E \lsem \bfZ \rsem \to  E[ \eps]/(\eps^2)$$ be the algebra  homomorphism
whose restriction to  $\bfZ$  is given by
\[
\kappa_{1+\eps}(u,z) = \chi_p(z)^{-1}(1 + \log_p( \chi_p(z)) \varepsilon),
\]
and let $E_{1+\eps}^{(p)}(\eta, \phi)$ be the image of $\cE(\eta, \phi)$  under $\kappa_{1+\eps}$. 
The Fourier expansion of $E_{1+\eps}^{(p)}(\eta, \phi)$ can be written as 
\begin{equation}
\label{eqn:Eis-fam}
E_{1+\eps}^{(p)}(\eta,\phi) \ = \ a_0(\eta,\phi) \ + \ \sum_{\nu}a_\nu(\eta,\phi) q^\nu.
\end{equation}
where the coefficients can  be  read off from \cite[Prop.~2.1, 3.2]{DDP11}, as summarised in the following lemma:
 \begin{lemma} 
\label{lemma:eisenstein-family}
The Fourier coefficients of $E_{1+\eps}^{(p)}(\eta,\phi)$ are given by
\begin{equation}
\label{eqn:Fourier_E_eta_psi}
\begin{array}{lll}
a_{0}(\eta,\phi) 
&=&  \ \displaystyle \frac{L_p'(\eta^{-1}\phi,0)}{4\eta(\fd)} \cdot \varepsilon 
\\ [7pt]
a_{\nu}(\eta,\phi)
&=&
\displaystyle 
\sum_{p \nmid I \mid (\nu) \fd} \eta\left( \frac{(\nu)\fd}{I}\right)\phi(I) \left(1 + \varepsilon \log_p \Nm(I) \right). 
\\ 
\end{array}
\end{equation}
\end{lemma}
 The article \cite{DDP11} constructs an explicit  {\em cuspidal}
family parametrised by $\cO_E\lsem \bfZ \rsem$ specialising to $f$ at weight 1. The family is not an eigenform over $\cO_E\lsem \bfZ \rsem$. Nevertheless, since $f$ itself is an eigenform,  one can deduce the existence of a morphism 
\begin{equation} 
 \label{eqn:def-pif}
\pi_f \colon \cT^{\nord} \to \cO_E 
\end{equation}
encoding the eigenvalues of Hecke operators acting on $f$. The composition with the morphism $\Lambda \to \cT^{\nord}$ will be denoted by 
\begin{equation}
\label {eqn: def-pi1}
\pi_1 \colon \Lambda \to \cO_E
\end{equation}
and corresponds to the character $\kappa_{\bft, \mathbf 0}\cdot \psi$. 

The remainder of \S~\ref{sec:deformation} will be dedicated to studying  lifts of the morphism $\pi_f$ to $E[\eps]/(\eps^2)$. Geometrically, this corresponds to studying the geometry of $\mathrm{Spec}(\cT^\nord)$ in an infinitesimal neighborhood of the prime ideal defined by $\pi_f$. 

\begin{remark}
The cuspidal family appearing in \cite{DDP11} was used to obtain the explicit formula for the derivative of the $p$-adic $L$-function $L_p(\psi, s)$ at $s=0$ conjectured by Gross, asserting that
 \begin{equation}
\label{eqn:Gross--Stark}
L_p'(\psi,0) = \sbL(\psi) L(\psi,0),
\end{equation}
where $\sbL(\psi)$ is the $\sbL$-invariant described in \S~\ref{subsec:Gal_prelim}.
In recent work of Betina, Dimitrov and Shih \cite {BDS} Gross' formula is linked to the study of the geometry of eigenvarieties from a Galois theoretic perspective. The  approach of \cite{BDS} informs the present work, and is carried out in a  broader setting.
\end{remark}

\subsection{Galois cohomology and $\sL$-invariants}
\label{subsec:Gal_prelim}

This section develops some results on Galois cohomology, which will be used 
later, notably  
 in \S~\ref{subsec:deformation_ring},
 to describe the tangent space of certain Galois deformation functors. These preliminary results are well known to experts. 
 
 Although it will not be used, it is worth noting that most of the arguments below are quite general and 
 also work for general number fields $F$. 
Let $H$ be a Galois extension of $F$ with Galois group $G$.
The $E[G]$-module $\Hom(G_{H}^\ab, E)$ can be described  explicitly via class field theory. This will be used to show that  certain global Galois cohomology classes for the totally odd character $\psi$
are determined by their images in  local cohomology. \\
A finite place $v$ of $H$ determines a prime ideal of $H$ whose decomposition group is isomorphic to the absolute Galois group of $H_v$, denoted by $G_v$.  For each $v$, there is an isomorphism between the completion of $H_v^\times$ and $G_{v}^\ab$ induced by the local Artin reciprocity map
\[
\mathrm {rec}_{v} \colon H_v^\times \ \lra \  G_{v}^\mathrm{ab}
\]
for which   the geometric normalisation is adopted. Note that the image of $G_v$ in $G_{H}^{\ab}$ is canonical. Thus, there is a restriction map $\res_v$   defined by
\begin{equation}
\label{lcft}
\begin{array}{lcll}
\res _v \colon & \Hom(G_H^\ab, E) & \lra & \Hom(H_v^\times, E)\\
	& f & \longmapsto & \mathrm {rec}_v \circ f|_{G_v^\mathrm{ab}}.
\end{array}
\end{equation}

The following lemma is well-known to experts, but   its proof is sketched below
for completeness. 
\begin{lemma}\label{l:exact sequence G_H^ab}
	There is an exact sequence of $E[G]$-modules
	\begin{equation}
		\label{e:exact sequence places|p}
		0 \to \Hom(G_H^\ab, E)\xrightarrow {(\res_v)_{v \mid p}}\prod_{v\mid p} \Hom(H_v^\times, E) \, \lra  \, \Hom(\cO_{H}[1/p]^\times, E).
	\end{equation}
	In addition, the rightmost map is surjective if and only if Leopoldt's Conjecture holds for $H$.

\begin{proof}
Let $\mathbb A_H^\times/H^\times$ be the id\`ele class group of $H$. The global Artin reciprocity map $\mathrm{rec}_H \colon \mathbb A_H^{\times} /H^{\times}  \to G_H^\ab$ is compatible with its local versions via the restriction maps, and gives a sequence 
\begin{equation}\label{e:exact sequence ideles}
0 \to \Hom(G_H^\ab, E)\to \Hom(\mathbb A_H^{\times}, E)\to \Hom(H^\times, E)
\end{equation}
 of continuous group homomorphisms and the topology on $H^\times$ is discrete. This sequence must be exact for topological reasons, and the two terms on the right are understood via their restrictions to the places above $p$, by the commutative diagram
\begin{center}
\begin{tikzcd}
\Hom(\mathbb A_H^\times, E) \arrow[r, "\Delta^\vee"] \arrow[d]
& \Hom(H^\times, E) \arrow[d] \\
\displaystyle \prod_{v\mid p} \Hom(H_v^\times, E) \arrow[r, "\Delta_p^\vee"]
& \Hom(\cO_H[1/p]^\times, E)
\end{tikzcd}
\end{center}
where $\Delta^\vee$, $\Delta_p^\vee$ are induced by the diagonal embeddings. Note that any element in $\Hom(\mathbb A_H^\times, E)$ must be trivial on the units $\cO_v^\times$ of $H_v$ for any $v\nmid p\infty$, and standard continuity arguments then show that the resulting map $\ker \Delta^\vee \to \ker \Delta_p^\vee$ is an isomorphism. Finally, by Dirichlet's Unit Theorem,
\[
\mathrm{rk}_\Z \, \cO_H[1/p]^\times=\mathrm{rk}_{\mathbb Z} \, \cO_H^\times+|\{v \mid p\}|,
\]
so the rightmost map of \eqref{e:exact sequence places|p} is surjective if and only if the $\Z_p$-rank of the closure of the image of $\cO_H^\times$ in $\prod_{v \mid p}\cO_v^\times$ is equal to the $\Z$-rank of $\cO_H^\times$, that is, if Leopoldt's Conjecture holds for $H$.
\end{proof}
\end{lemma}


Let $\phi \colon G \to E^\times$ be any character, viewed as a character of $G_F$, and consider the global
Galois cohomology group $\rH^1(F, E(\varphi)).$
The inflation-restriction sequence for continuous group cohomology of $E(\phi)$ leads to the identification
\begin{equation}\label{eqn:inflation-restriction}
\rH^1(F, E(\phi))\simeq \rH^1(H, E)^{\phi^{-1}} = \Hom(G_H^\mathrm{ab}, E)^{\phi^{-1}}.
\end{equation}
 The cohomology group $H^1(G, E(\phi))$ is related to certain units in $H$.  Let $S$ be a finite set of places of $F$, containing all infinite places and let $\cO_{H,S}^\times$ be the $S$-units of $H$, and write
\[
U_{\phi}:=(\cO_{H, S}^\times \otimes E)^{\phi^{-1}}.
\]
Then, the Galois-equivariant version of Dirichlet's Unit Theorem yields
\begin{equation}
	\label{Galois equivariant Dirichlet}
	\dim_E U_{\phi}=|\{ w \in S \mid \phi(w)=1 \}|-\dim_E E(\phi)^G.
\end{equation}

\medskip

The above discussion will now be specialised to the case where
 $H$ is the narrow Hilbert class field of $F$. In what follows, the character $\phi$ will be taken to be either the trivial character or   the unramified totally odd character $\psi$,
viewed as a character of $G=\Gal(H/F)$.  
 Let $S$ be the set of places containing the place corresponding to the prime $(p)$ of $F$ and all infinite places of $F$.
 It follows from 
\eqref{Galois equivariant Dirichlet}  that $\dim_E U_{\mathbbm 1}=d$, since $F$ is totally real and $\dim_E U_{\psi}=1$, because $(p)$ splits completely in $H/F$ and $\psi$ is totally odd.

Recall that we wrote
\[
 \alpha_1,\ldots, \alpha_d: F \hookrightarrow \bar\Q_p
\]
for the distinct $p$-adic embeddings of $F$. The prime ideal $p\cO_F$ splits completely in $H/F$, and the choice of a prime $\fp$ of $H$ above $p$ determines an identification $H_\fp = F_p$. Fix the choice of $\fp$ corresponding to the chosen embedding $\bar\Q \hookrightarrow \bar\Q_p$ once and for all, and write
\[
\tilde \alpha_1,\ldots,\tilde \alpha_d: H\hookrightarrow \bar\Q_p
\]
for the $p$-adic embeddings of $H$ extending $\alpha_1, \ldots, \alpha_d$ respectively, i.e.,
 \[
\tilde \alpha_{ j}|_{F}=\alpha_j \qquad \text{and} \qquad \tilde \alpha_{j}^{-1}(\fm _{\bar \Z_p})=\fp.
\]
 \begin{lemma}
\label{lemma:dim_coho} If the Leopoldt Conjecture holds for $H$, then: 
\begin{equation}
\dim_E \rH^1(F, E)=1 \qquad  \text{and} \qquad
\dim_E  \rH^1(F, E(\psi))=d.
\end{equation}
\begin{proof}
Combining \eqref{eqn:inflation-restriction}  with  Lemma \ref{l:exact sequence G_H^ab}, for any character $\phi$ of $G$,  the cohomology group $\rH^1(F, E(\phi))$ is equal to the $\phi^{-1}$-eigenspace of the kernel of
	\begin{equation} \label{e: Hom G_H^ab in TR_p case}
		\prod_{\fp_i\mid p} \Hom(H_{\fp_i}^\times, E) \lra  \Hom(\cO_{H}[1/p]^\times, E).
	\end{equation}
Let $\fp_1=\fp$ be the chosen prime. Then $\Hom(H_{\fp_1}^\times, E)$ is spanned by $(\log_p \circ \tilde  \alpha_j)_j$ and $\ord_{\fp_1}$. Since $p$ splits completely in $H/F$, the $E$-vector space
$\Hom(H_{\fp_1}^\times, E)$ generates the source of \eqref{e: Hom G_H^ab in TR_p case}  as an $E[G]$-module. The $E[G]$-span of each basis element of $\Hom (H_{\fp_1}^\times , E)$ is isomorphic to the right regular representation of $G$. In particular,  the multiplicity of every one-dimensional representation  of $G$ in the source of \eqref{e: Hom G_H^ab in TR_p case} is equal to $\dim_E (\Hom (H_{\fp_1}^\times , E))=(d+1)$.

Thus,  \eqref{Galois equivariant Dirichlet} leads to  the inequality (which is an equality if Leopoldt's Conjecture holds for $H$)
	\[
	\dim_E \rH^1(F, E(\phi)) \geq (d+1)-\dim_E U_{\phi^{-1}} =\begin{cases}
	1, \quad \text{if $\phi=\mathbbm 1,$}\\
	d, \quad  \text{if $\phi$ is totally odd}.
	\end{cases}
	\]
\end{proof}
\end{lemma}

\par When $\phi$ is either totally odd or trivial, we wish to describe the restriction to the decomposition group at the prime $p$ of a basis of $\rH^1(F, E(\phi))$. Denote
\[
\res_p \colon \rH^1(F, E(\phi)) \to \rH^1(F_p, E(\phi))
\]
for the restriction to the decomposition group at $p$, which is characterised by the choice of embedding $\bar \Q \hookrightarrow \bar \Q_p$. Since $\phi|_{G_{F_p}}=1$, by local class field theory, we can choose a basis of the target given by the $p$-adic valuation in $F_p^\times$,  denoted by $o_p$, together with the homomorphisms
\[
\ell_{p,j}=\log_p \circ \alpha_j, \qquad \qquad 1 \leq j \leq d.
\]

\par   Choose an auxiliary   prime $\mathfrak p_i$, for some $1 \leq i \leq n$. Then, under the running assumptions, $U_{\psi}$ is a one-dimensional $E$-vector space. Choose a generator $u_{\psi}$ of $U_{\psi}$, and note that $\ord_{\fp_i}(u_{\psi})\neq 0$, since the $\psi^{-1}$-eigenspace of  $\cO_H^\times \otimes E$ is trivial by \eqref{Galois equivariant Dirichlet}.
\begin{definition}
\label{def:L-invariant}
	The quantity
	\begin{equation}
		\mathscr L_j(\psi)=-\frac{(\log_p \circ \, \tilde\alpha_{j})(u_{\psi})}{(\ord_{p}\circ\tilde\alpha_{j})(u_{\psi})}, \qquad (1 \leq j \leq d)
	\end{equation}
is called the {\em partial $\sL$-invariant} of $\psi$ with respect to the $j$-th embedding  of $F$ into $\bar \Q_p$.
\end{definition}

Note that this definition is independent of the choice of the generator $u_{\psi}$ of $U_\psi$ and of the auxiliary choice of the prime $\fp$ of  $H$ above $p$. 
However, it depends on  the choice of the $p$-adic embedding $\alpha_{j}$ of $F$, thus justifying the notation. 

The following lemma related the partial $\sL$-invariants of Definition \ref{def:L-invariant} with the Gross--Stark
 unit $u_\tau$ of the introduction. In order to state it precisely, fix a choice of the unit  $u_\psi$ by fixing an RM point $\tau$ of discriminant $D$ and setting 
\begin{equation}
\label{eqn:def-u-phi}
u_\psi = \prod_{\sigma\in {\rm Gal}(H/F)} (\sigma u_\tau)^{\psi(\sigma^{-1})}.
\end{equation}
\begin{lemma}
\label{lemma:from-L-to-u}
For all odd characters $\psi$ of ${\rm Gal}(H/F)$, and all
$1\le j \le d$,
$$ \sL_j(\psi) L(\psi,0) = \log_p(\tilde\alpha_j(u_\psi)).$$
\begin{proof}
This follows after noting that, by definition of the Gross--Stark unit $u_\psi$, 
$$ \ord_p(\tilde\alpha_j(u_\psi)) = - L(\psi,0).$$
\end{proof}
\end{lemma}

The {\em full $\sbL$-invariant} of $\psi$ is the quantity
\[
\sbL(\psi)=\sum_{j=1}^d \sL_j(u_{\psi}). 
\]
It can alternatively be defined as
\[
\sbL(\psi)=-\frac{(\log_p \circ \Nm^{F_p}_{\Q_p}\circ \, \tilde\alpha_{j})(u_{\psi})}{
(\ord_p \circ \tilde\alpha_{j})(u_{\psi})}.
\]
 
\begin{remark}
$  $Of primary   interest is the case where $F$ is quadratic and $\psi$ is an odd narrow class character. Under theses assumptions, the character $\psi$ satisfies $\psi^{-1}(\sigma)=\psi(\tau \sigma \tau^{-1})$ for any $\tau \in G_{\Q}\smallsetminus G_{F}$.  This implies that the  partial $\sL$-invariants satisfy the relations
\[
\sL_1(\psi)=\sL_2(\psi^{-1}), \qquad  \sL_2(\psi) = \sL_1(\psi^{-1}).
\]

\end{remark}
\begin{prop}
	\label{p:dim cohomology groups}
	Let $\phi$ be any character of $G$. Denote $\chi_p \colon G_F \to \Z_p^\times$ the $p$-adic cyclotomic character. If Leopoldt's Conjecture holds for $H$, then:
	\begin{enumerate}
		\item If $\phi$ is trivial,  $\eta_{\mathbbm 1} := -\log_p \circ \chi_p$
		 generates $\rH^1(F, E(\phi))$.
		 Its restriction to the decomposition group at $p$ satisfies
		\[
		\res_p(\eta_{\mathbbm 1}) = \sum_{j=1}^d \ell_{p,j};
		\]
	\item If $\phi$ is totally odd, the  cohomology group $\rH^1(F, E(\phi))$ has a basis $\{ \eta_{\phi, j}\}_{1 \leq j \leq d }$ such that
		\[
		\res_p (\eta_{\phi, j}) = \ell_{p,j} + \mathscr L_j(\phi^{-1}) o_p.
		\]
 	\end{enumerate}
 \begin{proof}
By Proposition \ref{p:dim cohomology groups}, $\rH^1(F, E)$ is one dimensional; thus, it is generated by the (non-zero element) class of $\eta_{\mathbbm 1}=-\log_p(\chi_p)$. The restriction to the decomposition group at $p$ can be calculated by observing that, since $\chi_p$ is obtained by restriction to $G_F$ of a character of $G_{\Q}$, the same applies to the local characters at $p$. This implies that $\res_p (\eta_{\mathbbm 1})$ factors through the norm map from $F_p$ to $\Q_p$.

\par Since $\fp_1$ is the prime of $H$ determined by the fixed embedding $ \bar \Q \hookrightarrow \bar \Q_p$, the diagram

\begin{center}
    \begin{tikzcd}
		\rH^1(F, E(\phi))\arrow[r, "\res_p"] \arrow[d,"\res_H"]
		& \rH^1(F_p, E(\phi)) \arrow[d, "\res_{H_{\fp_1}}"] \\
		\rH^1(H, E(\phi)) \arrow[r, "\res_{\fp_1}"]
		&  \rH^1(H_{\fp_1}, E(\phi))
    \end{tikzcd}
\end{center}
commutes, where all the maps are given by restriction. In addition, $\res_{H_{\fp_1}}$ is an isomorphism, because $(p)$ splits completely in $H/F$; more precisely $\res_{H_{\fp_1}}$ satisfies
\begin{equation}
\label{e: from p to p_i}
\res_{H_{\fp_1}} (o_p)=\ord_p\circ \tilde\alpha_{j} \quad \text{and} \quad \res_{H_{\fp_1}} (\ell_{j,p})=\log_p\circ \tilde\alpha_{j}
\end{equation}
for every $ 1\leq j \leq d$.  It is worth noting at this stage
 that $\ord_p \circ \tilde\alpha_{j}$ is independent of $j$, while 
 $\log_p \circ \tilde\alpha_j$ depends very much on $j$.
After identifying $\rH^1(F, E(\phi))$ with the $\phi^{-1}$ eigenspace of the  kernel of \eqref{e: Hom G_H^ab in TR_p case}, the image of $\res_p$ is isomorphic to the image of this subspace via the projection
\[
 \bigoplus_{i=1}^n \Hom(H^\times_{\fp_i}, E) \lra \Hom(H^\times_{\fp_1}, E)
\]
on the first component
(which is, of course, not Galois equivariant). Let $\{1, 2, \dots, n\}$ be the $G$-set characterised by $\sigma \fp_i=\fp_{\sigma i}$ for every $\sigma \in G$.  Let $(f_i)_{i \in I}$ an element of $\bigoplus_{i=1}^n \Hom(H_{\fp_i}^\times, E)$. The action of $\sigma\in G$ on $(f_i)_{i}$ is given by
\[
\sigma(f_{i})_{i}=(f_{\sigma^{-1}i}\circ \sigma^{-1})_{i}.
\]
In particular, let $(f_i)_i$ belong to the $\phi^{-1}$-eigenspace. The action of $G$ on primes above $p$ is simply transitive. Let $i=\sigma^{-1}1$ for $\sigma \in G$. Then
\[
f_i=\phi(\sigma)^{-1} (f_1 \circ \sigma)
\]
Let $u_{\phi^{-1}}$ be a generator of the $\phi$-eigenspace of $\cO_H^\times \otimes E$. Then
\[
(f_i)_{i\in I}(\Delta_p(u_{\phi^{-1}}))=\sum_{i=1}^n f_i(u_{\phi})=\sum_{i=1}^n\phi(\sigma)^{-1}(f_1(\sigma u_{\phi}))
=nf_1(u_{\phi}).\]
Thus, $(f_i)_{i}$ belongs to the $\phi^{-1}$-component of the kernel of \eqref{e: Hom G_H^ab in TR_p case} if and only if $f_1(u_{\phi^{-1}})=0$. Let \[f_1=\sum_{j=1}^d x_j \log_p \circ \tilde\alpha_{ j}+y \ord_p\circ \tilde\alpha_{j}.
\]
The condition $f_1(u_{\phi^{-1}})=0$ cuts out a $d$-dimensional subspace, since $\ord_p\circ \tilde\alpha_{j}(u_{\phi^{-1}})\neq 0$. After re-writing this condition in terms of the $\sL$-invariants $\{\mathscr L_j(\phi^{-1})\}$ and comparing with \eqref{e: from p to p_i}, the proposition follows.
\end{proof}
\end{prop}

\subsection{$\Lambda$-adic Galois representations}
\label{subsec:HMF_prelim}

A general result of Hida establishes the connection between the nearly  ordinary Hecke algebra introduced in \S \ref{subsec:setup}  and Galois representations. More precisely, in \cite{Hid89}, certain  Galois representations are constructed which interpolate the representation corresponding to classical specialisations of Hida families for the Hecke algebra $\cT^{\nord}$. Exploiting the properties of these Galois representations, the study of the Hecke algebra $\cT^{\nord}$ infinitesimally at the prime ideal corresponding to the system of eigenvalues of $f$ can be reduced to that of a deformation ring that will be introduced in \S \ref{subsec:deformation_ring}.

The ultimate goal is to leverage the properties of the Galois representation to extract explicit formulae for the derivatives of the cuspidal family specialising to $f$, in the spirit of \cite{darmon-lauder-rotger1}. For this purpose, it suffices to consider the completed local ring $\cT_f$ obtained as the nilreduction of the completion of the localisation of $\cT^{\nord}$ at the prime ideal $\mathfrak q_f$ given by the kernel of the morphism $\pi_f$ defined in   \eqref{eqn:def-pif}. 
(Although this is not crucial for this application, it can be showed as in \cite[Prop. 6.4]{chenevier} that the completion of the localisation of $\cT^{\nord}$ at the point corresponding to $\mathfrak q_f$ is automatically reduced.)

It is natural to view $\cT_f$ as an algebra over $\Lambda_1$, the completion of the localisation of $\Lambda$ at the prime ideal $\mathfrak p_1=\ker \pi_1$; the latter is isomorphic to a ring of power series in $d+1$ variables over $E$.  

In this section, Hida's results are slightly refined in order to obtain a two-dimensional representation with coefficients in $\cT_f$ (Prop. \ref{p:existence galois stable lattice}), satisfying certain additional properties. The proof follows Mazur and Wiles' approach to the (somewhat delicate) study of deformations of residually reducible representations. The treatment of Bellaiche and Chenevier \cite{bellaiche-chenevier}, which is well-suited  to working over reduced henselian local rings such as $\cT_f$, will be followed. 

Write $K_f$ for the total ring of fractions of $\cT_f$; it is isomorphic to a product of fields $\prod_i K_{f,i}$, each corresponding to a minimal local component at of $\mathrm{Spec}(\cT^{\nord})$ at the point corresponding to $\mathfrak q_f$. 
 
Thus  $\cT_f$ can be viewed as a subring of $K_f$. In this context, the main result of \cite{Hid89} can be phrased as follows.

  \begin{theorem}[Hida] \label{Hida galois rep}
There exists a totally odd, continuous Galois representation 
\[
\rho_{K_f} \colon G_F  \to \GL_2(K_f)
\]
satisfying the following properties: 
	\begin{enumerate} 
	\item the pushforward $\rho_{K_{f_i}}$ of $\rho_{K_f}$  to $K_{f,i}$ is absolutely irreducible, for every i;
	\item $\rho_{K_f}$ is unramified outside $p$; 
	\item For every $\mathfrak l$ prime ideal of $F$ such that $\mathfrak l \nmid p$, let $\Frob_{\mathfrak l}$ be a Frobenius element. Then 
	\[
	\det (1-X \rho(\Frob_{\mathfrak l}))=1-\mathbf T_\mathfrak l X +\langle \mathfrak l \rangle \mathrm{Nm}(\mathfrak l) X^2;
	\]
	\item The restriction of $\rho_{K_f}$ to $G_{F_p}$ 
	is {\em nearly ordinary}, i.e., it satisfies
	\[
\rho_{K_f}|_{G_{F_p}} \simeq 
\left[
\begin{matrix}
\epsilon & * \\
0      & \delta 
\end{matrix}
\right],
	\]
	where 
	\begin{align}\label{e: local character}
	\delta \circ \rec_p(p) &= \mathbf{T}(p), 
	\\
	\delta \circ \rec_p(u) &=(\kappa^{\univ}) ((u,1)) \quad \mbox{for all } u \in \bfU,   
	\end{align}
where $\rec_p \colon F_p^\times \to G_{F_p}$ is the local Artin reciprocity map. 
\end{enumerate}	
\end{theorem}

 In order to relate the ring $\cT_f$ to a deformation ring, it is important to refine the Galois 
 representation $\rho_{K_f}$  to an integral version with coefficients in $\cT_f$. By the \v{C}ebotarev density theorem, the trace of $\rho_{K_f}$, as well as the characters $\delta$ and $\epsilon$, take values in $\cT_f\subset K_f$. Following Bella\"iche--Chenevier \cite{bellaiche-chenevier}, the existence of a free rank two $\cT_f$-module stable under the action of $G_F$ can be related to a certain global cohomology group. In addition, the  condition that the Galois representation $\rho_{K_f}$ is nearly ordinary imposes some local conditions on the global cohomology classes, that allow to show that  $\rho_{K_f}$ is conjugate to a representation with coefficients in $\cT_f$, following an argument which will now be described.
 
\par  
Let $M_{K_f}\simeq K_f^2$  be the two-dimensional Galois representation of $G_F$ provided 
by Theorem \ref{Hida galois rep}. Relative to a basis $(e^+,e^-)$  of $M_{K_f}$ on which a choice of complex conjugation for $F$ acts diagonally as
$\left [
\begin{smallmatrix}
 1 & 0\\
0 & -1
\end{smallmatrix}
\right ]$, 
the representation $\rho_{K_f}$ is given by                                                         
\[
\rho_{K_f}=\left[
\begin{matrix}
a_f & b_f\\
c_f & d_f
\end{matrix}                                                                                                                                                                                                                                           \right]
\colon G_F \to   \GL_2(K_f). 
\]
The fact that the traces of $\rho_{K_f}$ lie in $\cT_f$ implies that $a_f(\sigma) \pm d_f(\sigma)$ belong to $\cT_f$, and hence, that 
\[
a_f(\sigma), d_f(\sigma) \in \cT_f, \quad \mbox{ for all } \sigma\in G_F.
\]
It follows that
\[
\quad b_f(\sigma) \cdot c_f(\tau)  = a_f(\sigma\tau) - a_f(\sigma) a_f(\tau)   \in   \cT_f, \quad \mbox{ for all } \sigma, \tau\in G_F.
\]
Let $B_f$ and $C_f$ be the $\cT_f$-submodules of $K_f$ generated by the values of the functions $b_f$ and $c_f$ respectively. The {\em reducibility ideal}  $I^{\red}_f$  is the (proper) integral ideal of $\cT_f$ generated by the products $b_f(\sigma) c_f(\tau)$ for all $\sigma, \tau \in G_F$.

\par Fix a $G_{F_p}$-stable free one-dimensional submodule $L_{K_f}$ of $M_{K_f} = K_f^2$, and denote by $\epsilon_f$ and $\delta_f$ the local characters of $G_{F_p}$ acting on $L_{K_f}$ and $M_{K_f}/L_{K_f}$ respectively. 

\par
\begin{proposition}                                                                                                                                                                 \label{p:existence galois stable lattice}
There exists a free  $\cT_f[G_F]$-submodule $M_{\cT_f}\subset M_{K_f}$ of rank two over $\cT_f$, whose
 associated Galois representation $\rho_{\cT_f} \colon G_F \to \GL(M_{\cT_f})$ satisfies the following properties:
\begin{enumerate}[ wide = 0pt, leftmargin = *]
    \item The residual representation $M_{E}:=M_{\cT_f} \otimes E$ is semisimple;
    \item  \label{stable sub}
    There exists a free rank one summand $L_{\cT_f}$ of $M_{\cT_f}$ such that \begin{itemize}                                                                                  \item $L_{\cT_f}$ is $G_{F_p}$-stable and $G_{F_p}$ acts on $M_{\cT_f}/L_{\cT_f}$ via $\delta_f$;
    \item The subspace $L_{E}:=L_{\cT_f} \otimes E$ of $M_E$ is \emph{not} $G_{F}$-stable. 
    \end{itemize}
\end{enumerate}
\begin{proof}
	Let $B_{p,f}\subset B_f$   be the $\cT_f$-module generated by $b(G_{F_p})$, and likewise for $C_{f,p}\subset C_f$.
We claim that the natural   inclusion $ B_{p,f} \hookrightarrow B_f$ is
surjective.   The $\cT_f$-module $B_f$ is finitely generated by continuity of the representation $\rho_{K_f}$, and hence, by Nakayama's lemma, it suffices to show that the
induced  map
$$ i_B: B_{f,p}/\fm_{\cT_f}B_{f,p}   \rightarrow B_f/\fm_{\cT_f}B_f$$
is  surjective.
Consider the commutative diagram 
\begin{equation}\label{diagram B-B_p}
	\xymatrix{
		\Hom(B_f/\fm_{\cT_f} B_f, E)\ar[d]^{i_B^\vee} \ar[r]^{\qquad\Gamma} & \rH^1(F, E(\psi)) \ar[d]^{\res_p}\\	\Hom(B_{p,f}/\fm_{\cT_f} B_{p.f}, E)\ar[r]^{\quad \qquad  \Gamma_p}        & \rH^1(F_p, E)},
	\end{equation}
	where the top horizontal map $\Gamma$
	maps $\theta\in \Hom(B_f/\fm_{\cT_f}B_f, E)$ to the class of the cocycle
	\[
	\sigma\  \longmapsto \ \theta(b(\sigma))
	\]
	for every $\sigma \in G_F$, and $\Gamma_p$ is the corresponding map on local cohomology.  
	 By \cite[Lemma 3]{bellaiche-chenevier},
	the map $\Gamma$ is injective, and Proposition \ref{p:dim cohomology groups} implies that $\res_p$ is also injective. The commutativity of diagram \eqref{diagram B-B_p} implies that $i_B^\vee$ is injective, and therefore $i_B$ is surjective, as claimed, so that
	$B_f=B_{p,f}$. The same argument shows that 
	 $C_f=C_{p,f}$. 
	 
From the fact that $\rho_{K_{f_i}}$ is absolutely irreducible for every $i$, it follows that the the annihilator of the module $B_f$ (respectively $C_f$) is 0. 
One can deduce  that, without loss of generality, the vectors $(e^{+}, e^{-})$ can be rescaled by a pair of elements of $K_f^\times$ so that 
\[
L_{K_f}=\langle e^{+}+e^{-} \rangle.
\] 
Note that this basis is unique up to scaling, and hence, the resulting matrix representation of 
$ \rho_{K_f}$ is uniquely determined.
	
Changing the basis $(e^+,e^-)$ to $( e^{-}+e^{+}, e^+)$, the representation $\rho_{K_f}$ is given  in matrix form 
by
\[
	\left [
	\begin{matrix}
	1 & 0\\
	-1 & 1
	\end{matrix}
	\right ]
	\left [
	\begin{matrix}
	a_f & b_f\\
	c_f & d_f
	\end{matrix}
	\right ]
	\left [
	\begin{matrix}
	1 & 0\\
	1 & 1
	\end{matrix}
	\right ]
	=
	\left [
	\begin{matrix}
	a_f+ b_f & b_f\\
	-b_f+(d_f-a_f)+ c_f & - b_f+d_f
	\end{matrix}
	\right ].
	\]
	In particular, for every $\sigma \in G_{F_p}$, 
	\begin{equation*}
    a_f(\sigma)+b_f(\sigma)=\epsilon_f(\sigma)\qquad \text{and} \qquad -b_f(\sigma)+d_f(\sigma)-a_f(\sigma)+ c_f(\sigma)=0. 
	\end{equation*}
	Since $a_f$ and $d_f$ are valued in $\cT_f$ and $a_f(\sigma)=d_f(\sigma)=1 \pmod {\fm_{\cT_f}}$, the first equation implies that $b_{f}$ takes values in $\fm_{\cT_f}$.  Similarly, because $d_f(\sigma)-1\in \fm_{\cT_f}$, from the second equation it also follows that $c_f(\sigma)\in \fm_{\cT_f}$. Thus, $B_{f}=B_{p,f} \subset \fm_{\cT_f}$ and  $C_{p}=C_{p,f} \subset \fm_{\cT_f}$.	It follows that $\rho_{K_f}$  has coefficients in $\cT_f$ relative to the basis $(e^+,e^-)$, giving rise to a Galois stable $\cT_f$-lattice  $M_{\cT_f}:= \cT_f e^+ + \cT_f e^-$  which satisfies all the claims of Proposition \ref{p:existence galois stable lattice}.
\end{proof}
\end{proposition}

\subsection{A deformation ring for residually reducible representations}
\label{subsec:deformation_ring}

This section describes an abstract deformation ring $\cR^\nord_{\rho, L}$ relevant to the Eisenstein series $f$ defined in \eqref{p-stabilisation}.
This deformation ring is equipped with
a natural   $\Lambda_{1}$-algebra structure
\begin{equation} 
\label{Lambda_1 structure on R}
 \Phi_{\cR} \colon \Lambda_{1} \ \lra \ \cR^\nord_{\rho, L}.
\end{equation}
 The construction 
 of  $\cR^\nord_{\rho, L}$ is complicated by the residual reducibility, and we follow the approach of Calegari--Emerton \cite{calegari-emerton} to overcome this. The main result of this section is Proposition \ref{l:explicit infinitesimal deformation}, which computes the map induced by $\Phi_{\cR}$ on tangent spaces and shows it is an isomorphism.

\subsubsection{The deformation functor. }
Let $\mathfrak C_E$ be the category of local complete noetherian rings with residue field $E$. Consider the functor 
\[
\cD^{\det} \times \cD^{\loc} \colon \mathfrak C_E \, \lra \, \mathbf {Sets}
\]
which takes any $(A, \mathfrak m_A) \in \mathrm{Ob}(\mathfrak C_E)$ to the pairs of continuous  characters
\[
(\upsilon_A,  \vartheta_A) \colon \bfU \times \bfZ  \to A^\times , \quad\text{ with } \quad (\upsilon_A, \vartheta_A)=(1, 1) \pmod {\fm_A}.\]

\par This functor is represented by the completed local ring $\Lambda_{1}$ of the Iwasawa algebra introduced in \S~\ref{subsec:HMF_prelim}. The ring $\Lambda_{1}$ is isomorphic to a ring of formal power series in $(d+1)$ variables over $E$.

 \par We first define a deformation ring $\cR_{\rho, L}$. Let $\rho=\psi \oplus 1$ with the standard basis $(v_1, v_2)$ of $V=E^2$. This is a semisimple reducible representation, and as such it admits non-scalar endomorphisms. To obtain a representable functor, the  deformation problem needs to be  suitably rigidified. This is done, following Calegari--Emerton \cite{calegari-emerton},   by setting $L=\langle v_1 + v_2 \rangle$, which is a line that is \textit{not} stable under the action of $G_F$. 
A {\em strict deformation} of $(V,L,\rho)$ 
over an $E$-algebra $A$  is a 
quadruple  $(V_A, L_A, \rho_A, g_A)$ where 
\begin{itemize}
\item $V_A$ is a free $A$-module of rank $2$;
\item $L_A$ is a free rank $1$ summand of $V_A$;
\item $\rho_A \colon G_{F} \to \GL(V_A)$ is a continuous representation;
\item $g_A \colon V_A \otimes_A E \simeq V$ is an isomorphism of $E[G_\Q]$-modules sending  $L_A \otimes E$ to $L$.
\end{itemize}
Two strict deformations $(V_A, L_A, \rho_A,g_A)$ 
and $(V_A', L_A', \rho_A',g_A')$ are said to be
{\em equivalent} if there is an $A[G_\Q]$-module 
isomorphism $h: V_A\lra V_A'$ sending $L_A$ to $L_A'$ and for which the diagram
$$
\xymatrix{
V_A  \otimes_A E  \ar@{->}[r]^{\quad g_A} \ar[d]_{h\otimes E} & V
\ar@{=}[d]  \\
V_A'  \otimes_A E  \ar@{->}[r]^{\quad g_A'} & V }
$$
commutes. A {\em deformation} of $(V,L,\rho)$ over an $E$-algebra $A$ is an equivalence class of strict deformations over $A$.

Consider the functor
\[
\cD_{\rho, L} \colon \mathfrak C_E \, \lra \, \mathbf{Sets}
\] 
which associates to an object $A$ in $\mathfrak C_E$ the set of deformations of $(V,L,\rho)$ over $A$.
The functor $\cD_{\rho, L}$ is representable by a complete local Noetherian ring  $\cR_{\rho, L}$ with residue field $E$. 
The representability can be verified as in \cite[Prop.~1]{mazur}. It is ensured by the additional  datum of a  line lifting $L$, 
which can be viewed a ``partial framing" of the functor parametrising deformations, forcing   automorphisms of a deformation  to consist only of scalars.  Thus rigidified, the  deformation functor presents the advantage of being fine enough to be representable, while still being conceivably comparable with a Hecke algebra, as in \cite{calegari-emerton}.

\par Finally,   consider the functor $\cD^\nord_{\rho, L}$ classifying quintuples $(V_A, L_A, \rho_A, g_A, \iota_A)$ such that: 
\begin{itemize}
\item the equivalence class of the quadruple $(V_A, L_A, \rho_A, g_A)$ belongs to $\cD_{\rho, L}(A)$;
\item the free rank one summand  $L_A$ is $G_{F_p}$-stable; 
\item $\iota_A \colon F_{p}^\times \to A^\times$ is the character satisfying 
\[
((\rho_A \circ \rec_p)y)v=\iota_A(y)v \mod L_A
\] 
for every $y \in F_p^\times$ and $v\in V_A$, where $\rec_p$ denotes the local Artin reciprocity map. 
\end{itemize}
\begin{remark}
Of course the datum of a character $\iota_A \colon F_{p}^\times \to A^\times$ is redundant in the previous definition; it is nonetheless useful to keep track of it, since it plays a key role in the calculations. 
\end{remark}

The deformation functor $\cD^\nord_{\rho, L}$ is representable by a quotient of $\cR_{\rho, L}$, denoted 
$\cR^\nord_{\rho, L}$. Indeed,   choose an 
$\cR_{\rho, L}$-basis $(\tilde v_{1}, \tilde v_{2})$ for the universal representation, lifting $(v_1,v_2)$, in such a way that the universal free rank one summand is given by
\[
L_{\cR_{\rho, L}} = \langle \tilde v_{1}+\tilde v_{2} \rangle.
\] 
Then the  ring  $\cR^\nord_{\rho, L}$ is the quotient of $\cR_{\rho, L}$ by the ideal
\[
J^\nord=\langle \alpha(\sigma)+\beta(\sigma)-\gamma(\sigma)-\delta(\sigma) \mid \forall \sigma \in G_{F_p}\rangle,
\]
where $\alpha,\beta,\gamma,\delta$ are the entries of $\rho_{\cR_{\rho, L}}=\left[
\begin{smallmatrix}
\alpha & \beta\\
\gamma & \delta
\end{smallmatrix}
\right]$ with respect to the chosen basis. 

\par We now describe the $\Lambda_{1}$-algebra structure of $\cR^\nord_{\rho, L}$. Given a quintuple $(V_A, L_A, \rho_A, g_A, \iota_A) \in \cD^\nord_{\rho, L}(A)$, the pair
$
(\det \rho_A \cdot \psi^{-1}, \iota|_{\cO_{F_p}^\times} )
$
belongs to $\cD^{\det} \times \cD^\loc(A)$. In particular,   taking (a representative of) the universal object of $\cD^\nord_{\rho, L}$ yields a morphism $\Phi_{\cR}$.
As a consequence of Proposition \ref{p:existence galois stable lattice}, the universal property of $\cR^\nord_{\rho, L}$ yields a morphism $\Upsilon$ from the deformation ring to the Hecke algebra, which makes the following diagram commute: 
\begin{equation}
\label{R=T diagram}
\xymatrix{
	\cR^{\nord}_{\rho, L} \ \ar[rr]^{\Upsilon}  &  &    \quad      {\cT_f} \\
 & \Lambda_{1} \ar[lu]^{\Phi_{\cR}} \ar[ru]_{\Phi_{\cT}} & }
\end{equation}

\subsubsection{Tangent spaces. }
We now come to the main results of this subsection, and describe the map on tangent spaces induced by
 $\Phi_R$. Let $E[\eps]$ be the ring of dual numbers over $E$. Let
\[
t^{\det, \loc} = (\cD^{\det} \times \cD^\loc) (E[\eps]), \qquad 
 t_{\rho, L}=\cD_{\rho, L}(E[\eps])  \qquad \text{and} \qquad t^\nord_{\rho, L}=\cD^\nord_{\rho, L}(E[\eps])
\] 
be the tangent spaces of the three functors introduced above. The following lemma describes them explicitly.

\begin{lemma}
	\label{tangent t^nord}
	\begin{enumerate}[ wide = 0pt, leftmargin = *]
	\item There is an
	isomorphism
		 $G:\rH^1(F, \ad(\rho)) \to t_{\rho,L}$ sending the cohomology class of $\left[
\begin{smallmatrix}
\alpha & \beta \\
\gamma & \delta
\end{smallmatrix} \right]$ in $\rH^1(F, \ad (\rho))$ to the equivalence class of
\begin{equation}\label{def G}
\left(E[\eps]^2,\left( 1+\eps
\left[ \begin{smallmatrix} \alpha & \beta \\ \gamma & \delta \end{smallmatrix} \right ]\right )\rho,
\left \langle \left[
\begin{smallmatrix}
1 \\1
\end{smallmatrix}
\right]
\right \rangle, g_{\eps} \right),
\end{equation}
where $g_{\eps}$ sends the standard basis of $E[\eps]^2$ to $(v_1, v_2)$.\\
 \item Let $\rH^1(F, \ad (\rho))^\nord$ be the subspace of $\rH^1(F, \ad (\rho))$ consisting of cocycles for which
\[
\alpha(\sigma)+\beta(\sigma)-\gamma(\sigma)-\delta(\sigma) = 0,  \qquad \forall \sigma \in G_{F_p}.\]
There is an isomorphism $G' \colon \rH^1(F, \ad (\rho))^{\nord} \to t^{\nord}_{\rho, L}$ given by 
\[
G' \left( \left[ \begin{smallmatrix} \alpha & \beta \\ \gamma & \delta \end{smallmatrix} \right ] \right )=\left (G\left( \left[ \begin{smallmatrix} \alpha & \beta \\ \gamma & \delta \end{smallmatrix} \right ] \right ), 1+\eps (\delta-\beta) \circ \rec_p )\right ).
\]
	\end{enumerate}
\begin{proof}
\begin{enumerate}[ wide = 0pt, leftmargin = *]
	\item
For any equivalence class in $t_{\rho, L}=\cD_{\rho, L}(E[\eps])$, we can choose a representative of the form \eqref{def G} for some cocycle in $\rZ^1(F, \ad (\rho))$. It suffices to verify that $G$ is well-defined; in other words, it is enough to show that two lifts if $\rho_{\eps}, \rho'_\eps$ of $\rho$ are conjugate by a matrix \[g\in \ker (\GL_2(E[\eps]) \to \GL_2(E)),\] then they are conjugate by a matrix stabilizing the line $\left \langle \left[
\begin{smallmatrix}
	1 \\1
\end{smallmatrix}
\right]
\right \rangle$.
This follows   from the fact that the space of coboundaries for the adjoint representation of the form
\[
\sigma \mapsto
\left[
\begin{matrix}
\psi(\sigma) & 0 \\
0 & 1
\end{matrix}
\right ]
 \left[
	\begin{matrix}
r & s \\
t & u
\end{matrix}
\right ]
\left[
\begin{matrix}
\psi(\sigma)^{-1} & 0 \\
0 & 1
\end{matrix}
\right ]- \left[
\begin{matrix}
r & s \\
t & u
\end{matrix}
\right ]=
 \left[
\begin{matrix}
0 & (\psi(\sigma)-1)s \\
(\psi(\sigma^{-1})-1)t & 0
\end{matrix}
\right ]
\]
for $\sigma \in G_F$ and $ \left[
\begin{smallmatrix}
r & s \\
t & u
\end{smallmatrix}
\right ] \in \mathrm{M}_2(E)$ is spanned by coboundaries of matrices fixing the line $\left \langle \left[
\begin{smallmatrix}
	1 \\1
\end{smallmatrix}
\right]
\right \rangle$.

\item
Let $X$ be the cohomology class of $\left[\begin{smallmatrix}
\alpha & \beta\\
\gamma & \delta
\end{smallmatrix}
\right] \in \rZ^1(F, \ad (\rho))$. Denote $\rho_{\eps}=(1+\eps \left[\begin{smallmatrix}
\alpha & \beta\\
\gamma & \delta
\end{smallmatrix}
\right])\rho$.
The image of $X$ under $G$ belongs to the $t_{\rho, L}^\nord$ if the line $\left \langle \left[
\begin{smallmatrix}
	1 \\1
\end{smallmatrix}
\right]
\right \rangle$ is stable for the action of $G_{F_p}$. In the basis $\left (\left \langle \left[
\begin{smallmatrix}
1 \\1
\end{smallmatrix}
\right]
\right \rangle, \left \langle \left[
\begin{smallmatrix}
	0\\1
\end{smallmatrix}
\right]
\right \rangle
\right )$
the matrix of $\rho_{\eps}$ is given by
\[
\left[\begin{matrix}
1 & 0\\
1 & 1
\end{matrix}
\right]^{-1}
\cdot \left[\begin{matrix}
\alpha & \beta\\
\gamma & \delta
\end{matrix}\right]
\cdot \left[\begin{matrix}
1 & 0\\
1 & 1
\end{matrix}
\right]
=
\left[\begin{matrix}
\psi+\eps(\psi \alpha+\beta) & \beta \eps \\
(1-\psi)+\eps(-\alpha \psi-\beta+\gamma \psi +\delta)	 & 1+\delta \eps-\beta \eps
\end{matrix}
\right].
\]
The claim follows by observing that the restriction of $\psi$ to $G_{F_p}$ is trivial.
\end{enumerate}
\end{proof}
\end{lemma}

\par With these identifications, explicit bases for the relevant tangent spaces can be obtained by using the results in \S~\ref{subsec:Gal_prelim}. Identify $t^{\det, \loc}$ with $\Hom(\bfZ, E) \oplus \Hom (\bfU, E)$,
and choose the basis 
\begin{equation}
e_0:=(\eta_{\mathbbm 1}, 0), \qquad e_{j}=(0, \ell_{p,j}), \qquad 1 \leq j \leq  d. 
\end{equation}

Given this choice of basis of the tangent space, an identification between $\Lambda_1 = E \lsem X_0, X_1, \dots, X_d \rsem$ can be chosen so that the universal characters are given by the pair $(\upsilon^\univ, \vartheta^\univ)$ satisfying
\begin{align*}
\vartheta^\univ =1+X_0e_0 \pmod {\fm^2_{\Lambda_1}} && \upsilon^\univ=1+ \sum_{j=1}^d X_je_j \pmod {\fm^2_{\Lambda_1}} .
\end{align*}

For the tangent space $t_{\rho,L}$, note  that since $\rho=\psi \oplus \mathbbm 1$,  it follows that
\[
\ad(\rho) = \mathbbm 1^2 \oplus \psi \oplus \psi^{-1}. 
\]
i.e. the adjoint representation of $\rho$ splits completely. Hence, by Lemma \ref{lemma:dim_coho},   the cohomology group $\rH^1(F, \ad (\rho))$ has dimension $2d+2$, and we may choose the $E$-basis consisting of
\begin{equation}
\label{e: basis H^1(F,ad rho)}
A=\left [\begin{matrix}
\eta_{\mathbbm 1} & 0\\
0 & 0
\end{matrix}
\right],
\quad
D=\left [\begin{matrix}
0 & 0\\
0 & \eta_{\mathbbm 1}
\end{matrix}
\right],
\quad
B_j=
\left [\begin{matrix}
0 & \eta_{\psi, j}\\
0 & 0
\end{matrix}
\right],
\quad
C_j=\left [\begin{matrix}
0 & 0 \\
\eta_{\psi^{-1},j} & 0
\end{matrix}
\right], \qquad 1 \leq j \leq d. 
\end{equation}
where the entries are described by Proposition \ref {p:dim cohomology groups}. With respect to these choices of bases, we now explicitly describe the map on tangent spaces induced by $\Phi_{\cR}$, denoted
\begin{equation} 
\label{e: tangent structural R}
\Theta  \colon t^\nord_{\rho, L} \, \lra \, t^{\det, \loc}. 
\end{equation}

\begin{proposition} 
\label{l:explicit infinitesimal deformation}
If $\sbL(\psi)+ \sbL(\psi^{-1}) \neq 0$, the map $\Theta$ is an isomorphism and its inverse satisfies
\begin{align}
\Theta^{-1}(e_0)&=
\frac{\mathscrbf L(\psi^{-1})}{\mathscrbf L(\psi)+\mathscrbf L(\psi^{-1})}
\left ( A+\textstyle \sum_{k=1}^d C_k \right )+
\frac{\mathscrbf L(\psi)}{\mathscrbf L(\psi)+\mathscrbf L(\psi^{-1})}
 \left (D+ \textstyle \sum_{k=1}^d B_k  \right )\\
\Theta^{-1}(e_j)&=
\frac{\mathscr L_j(\psi)-\mathscr L_j(\psi^{-1})}{\mathscrbf L(\psi)+\mathscrbf L(\psi^{-1})} (A-D+\textstyle \sum_{k=1}^d C_k-\sum_{k=1}^d B_k)-B_j-C_j
\end{align}
for $1 \leq j \leq d$.
\begin{proof}
By Lemma \ref{tangent t^nord}, up to composing with $G'$,    the tangent space $t^\nord_{\rho, L}$ is identified with a subspace of $\rH^1(F, \ad (\rho))$ given by the kernel of
\begin{align*}
P_1 \colon \rH^1(F, \ad (\rho) ) \to \rH^1(F_p,E), &&
\left [
\begin{smallmatrix}
\alpha & \beta \\
\gamma & \delta
\end{smallmatrix}
\right ]
\mapsto \res_{p}(\alpha+\beta-\gamma-\delta).
\end{align*}
where $\res_{p}$ denotes the restriction to $G_{F_p}$ (note that this is well defined, because $\psi(G_{F_p})=1$).
Similarly, the map $\Theta\colon t^\nord_{\rho, L} \to t^{\det, \loc}$
sends a quintuple $(V_\eps, L_\eps, \rho_{\eps}, g_{\eps},\iota_{\eps})$ to the pair $(\det \rho_{\eps}\cdot \psi^{-1}, \iota_{\eps}|_{\cO_{F_p}^\times })$. Again, by Lemma \ref{tangent t^nord}, it can be interpreted in terms of Galois cohomology as the restriction to the kernel of $P_1$ of
\begin{align*}
P_2 \colon \rH^1(F, \ad (\rho) ) \to \Hom(\bfZtor, E)\oplus \Hom(\bfUtor, E), && 
\left [
\begin{smallmatrix}
\alpha & \beta \\
\gamma & \delta
\end{smallmatrix}
\right ]  \mapsto \left(\alpha+\delta,(\delta-\beta)\circ \rec_{p}|_{\cO_{F_p}^\times} \right).
\end{align*}
Thus, in order to show that $\Theta$ is an isomorphism, it suffices to show that $P:=(P_1,P_2)$ is an isomorphism. Choose the bases $(A, D, B_1, \dots, B_d, C_1, \dots, C_d)$ for $\rH^1(F, \ad (\rho))$ and
\begin{align*}
(o_p,0), \ (\ell_{p, 1},0), \dots, (\ell_{p,d},0) , \ (0, e_0), \dots,(0,e_d)
\end{align*}
for the target of $P$.
Let $(x_{A}, x_D, x_{B_1}, \dots, x_{B_d}, x_{C_1}, \dots, x_{C_d})$ be the coordinates of a class in $\rH^1(F, \ad (\rho))$ with respect to the basis above. The map $P$ yields a system of $(2d+2)$ linear equations in $(2d+2)$-variables: 
\begin{align*}
\textstyle \sum_{j=1}^d (\sL_j(\psi)x_{C_j} -\sL_j(\psi^{-1})x_{B_j})=0\\
x_A-x_D+x_{B_j}-x_{C_j}=0\\
x_A+x_D=0\\
x_D-x_{B_j}=0
\end{align*}
where $1 \leq j \leq d$. The corresponding matrix has determinant $\sbL(\psi)+\sbL(\psi^{-1})$. 
The expressions for the inverse of $\Theta$ can be obtained by the inverting the matrix of $P$ with respect to the above  bases.
\end{proof}
\end{proposition}

\begin{remark}
In the case of interest to this paper where $[F:\Q] =2$, the non-vanishing
of $\sbL(\psi)+\sbL(\psi^{-1})$   is clear since $\sbL(\psi) = \sbL(\psi^{-1}) \neq 0$. In fact, this non-vanishing holds in general. For details, compare with \cite{BDS}.
\end{remark}
Recall the commutative
diagram 
\eqref{R=T diagram} arising from 
Proposition \ref{p:existence galois stable lattice}.

\begin{theorem} \label{thm: R=T}
The maps $\Phi_{\cR}, \Phi_{\cT}$ and $ \Upsilon$ are isomorphisms.  

\begin{proof}
Since $\Phi_{\cR}$ is a morphism of complete local noetherian rings with residue field $E$, the injectivity of \eqref{e: tangent structural R} implies that $\Phi_{\cR}$ is surjective. The top row is surjective because all Hecke operators are in the image; thus it follows that $\Phi_{\cT}$ is surjective. But $\cT_f$ is a torsion free $\Lambda_1$-algebra; in particular $\Phi_{\cT}$ is an injective, hence an isomorphism. It follows that $\Phi_{\cR}$ and $\Upsilon$ are isomorphisms as well.
\end{proof}
\end{theorem} 

\begin{prop}
\label{prop:lambda-mu-xi}
 The inverse of $\Phi_{\cT}$ satisfies 
 \begin{equation}
 \begin{array}{llll}
 \Phi_{\cT}^{-1}(\mathbf T_{\mathfrak l})&=&1+\psi(\mathfrak l)+\log_p(\mathrm{Nm}(\mathfrak l))\cdot (\lambda+\mu \psi(\mathfrak l)) & \pmod{ \mathfrak m_{\Lambda_{1}}^2 }\\
 \Phi_{\cT}^{-1}(\langle \mathfrak l \rangle) \Nm(\mathfrak l)&=& \psi(\mathfrak l)(1+(\lambda+\mu)\log_p(\mathrm{Nm}(\mathfrak l))) & \pmod {\fm_{\Lambda_1}^2 }\\
  \label{T_p}
 \Phi_{\cT}^{-1}(\mathbf T(p)) &=&1+\xi  & \pmod {\mathfrak m_{\Lambda_{1}}^2},
 \end{array}
\end{equation}
  where $\mathfrak l$ is a prime ideal of $F$ such that $p \nmid \rm{Nm}(\mathfrak l)$ and $\lambda, \mu, \xi \in \mathfrak m_{\Lambda_{1}}$ are given by 
\[
\begin{array}{lll}
\lambda &=& (\sbL(\psi)+\sbL (\psi^{-1}))^{-1} \left( \sbL(\psi^{-1})X_0+\sum_{j=1}^d(\sL_j(\psi)-\sL(\psi^{-1}))X_j\right)\\
 \mu &=& (\sbL(\psi)+\sbL (\psi^{-1}))^{-1} \left( \sbL(\psi)X_0-\sum_{j=1}^d(\sL_j(\psi)-\sL(\psi^{-1}))X_j\right)\\
\xi &=& (\sbL(\psi)+\sbL (\psi^{-1}))^{-1} \left( -\sbL(\psi)\sbL(\psi^{-1})X_0+\sum_{j=1}^d(\sL_j(\psi)\sbL(\psi^{-1})+\sL_j(\psi^{-1}) \sbL(\psi))X_j\right).
\end{array}
\]

\begin{proof}
Let $\left(V^{\univ}, L^{\univ}, \rho^{\univ}, g^\univ, \iota ^{\univ}\right)$ be a representative of the universal object of the functor $\cD^\nord_{\rho, L}$ over the deformation ring $\cR^\nord_{\rho, L}$. Then  
	\[
	\Upsilon^{-1}(\mathbf T_{\mathfrak l})=\Tr(\rho^{\univ}) (\Frob_{\mathfrak l}) \qquad \text{and} \qquad   \Upsilon^{-1}(\mathbf T(p))=\iota^{\univ}(\Frob_{p}).	
	\]
Denoting 
\[
\left (\begin{matrix}
	a & b \\
	c & d
	\end{matrix} \right ) = \Phi^{-1}_{\cR} \circ \rho^\univ \colon G_F \ \lra \ \GL_2(\Lambda_{1}),
\]
it follows  from Proposition \ref{l:explicit infinitesimal deformation} that modulo $\fm_{\Lambda_{1}}^2$,  
\[
\begin{array}{lll}
a &=& \psi+ \psi \eta_{\mathbbm 1}(\sbL(\psi)+\sbL (\psi^{-1}))^{-1} 
\left( \sbL(\psi^{-1})X_0+\sum_{j=1}^d(\sL_j(\psi)-\sL(\psi^{-1}))X_j\right)\\
d &=& 1+\eta_{\mathbbm {1}}(\sbL(\psi)+\sbL (\psi^{-1}))^{-1} 
\left( \sbL(\psi)X_0-\sum_{j=1}^d(\sL_j(\psi)-\sL(\psi^{-1}))X_j\right) 
\end{array}
\]
from which  the expression for $\Upsilon ^{-1}(\mathbf T_{\mathfrak l})$ is obtained. Similarly,   $\Phi_{\cR}^{-1} \circ \Phi^\univ=d-b$. Since   $\eta_{\mathbbm 1}(\Frob_p)=0$, it can be seen that $(d-b)(\Frob_p)=1-b(\Frob_p)$, which is equal to
\[
1+\left(\frac{\sbL(\psi)(\sum_{k=1}^d\eta_{\psi,k}(\Frob_p))}{\sbL(\psi)+\sbL(\psi^{-1})}
X_0+
\sum_{j=1}^d  {\frac{(\sL_j(\psi^{-1})-\sL_j(\psi))\sum_{k=1}^d\eta_{\psi,k}(\Frob_p)}{\sbL(\psi)+\sbL(\psi^{-1})}}-\eta_{\psi, j}(\Frob_p)
\right ) X_j
\]
modulo $\fm^2_{\Lambda_{1}}$ by Proposition \ref{l:explicit infinitesimal deformation}. The equality \eqref{T_p} then follows from Lemma \ref{p:dim cohomology groups}. 
\end{proof}
\end{prop}

\subsection{Fourier coefficients}
\label{subsec:q-exp}
The above results will now be specialised to the case where $F$ is a real quadratic field and $\psi$ is an unramified totally odd character of $F$. In this case, we compute the Fourier coefficients of the \textit{anti-parallel} family through the Eisenstein series $f$ of parallel weight 1  discussed above. 

\par The \textit{anti-parallel} weight direction in the tangent space of the Iwasawa algebra  is  the direction corresponding to the morphism $\Lambda_{1} \lra E [\varepsilon]/ (\varepsilon^2)$ given in terms of generators as follows:
\begin{equation}\label{eqn:antiparallel direction} 
\left\{ 
\begin{array}{cll}
X_0, X_1 & \longmapsto & \varepsilon, \\
X_2 & \longmapsto & 0.
\end{array}
\right.
\end{equation} 
Since, by Theorem \ref{thm: R=T}, the structural map $\Phi_{\cT} \colon\Lambda_1 \to \cT_f$ is an isomorphism, the map \eqref{eqn:antiparallel direction}  gives rise to a morphism from the nearly ordinary Hecke algebra to the ring of dual numbers
\begin{equation}
\tilde \pi_f : \cT^{\nord} \lra E[\varepsilon]/ (\varepsilon^2).
\end{equation}
lifting the morphism $\pi_f$ defined in \eqref{eqn:def-pif}.
This corresponds to the system of Hecke eigenvalues of a  first order eigenfamily of Hilbert modular forms $\cF$, whose Fourier coefficients, which can be recovered from $\tilde \pi_f$,  play a central role in what follows.

\par In the anti-parallel weight direction, one immediately checks that the quantities $\lambda,\mu$, and $\xi$ appearing in the description of the tangent space of $\cT^{\nord}$ in Proposition \ref{prop:lambda-mu-xi} specialise to
\begin{equation}
\lambda = \frac{\sL_1(\psi)}{\sbL(\psi)}  \qquad \qquad \mu =\frac{\sL_2(\psi)}{\sbL(\psi)} \qquad \qquad \xi = 0.
\end{equation}
Using the results in \S~\ref{subsec:deformation_ring}, the image under the anti-parallel weight morphism $\tilde \pi_f$ of the operators $\langle l \rangle$ and $\mathbf{T}_{\fl^n}$ for $\fl \neq (p)$, as well as $\mathbf{T}(p^n)$,
can now be  computed  in terms of these quantities:
\begin{itemize}
\item Let $\fl \neq (p)$ be a prime ideal.  
 Proposition \ref{prop:lambda-mu-xi} immediately implies  that
\[ 
\begin{array}{llcll}
\tilde \pi_f(\mathbf{T}_{\fl}) &=& 1 + \psi(\fl) & + & \varepsilon \log_p \left(\Nm(\fl)\right) \cdot \left(\lambda \psi(\fl) + \mu \right)\\
\tilde \pi_f(\langle \fl \rangle)\Nm(\fl) &=& \psi(\fl) & + & \varepsilon \log_p \left(\Nm(\fl) \right)\cdot \psi(\fl).
\end{array}
\]
The recursion relation proved in \cite[Corollary 4.2]{Hid88}, which states
\[
\mathbf{T}_{\fl^{n}} \mathbf{T}_{\fl} \ = \ \mathbf{T}_{\fl^{n+1}} \ + \ \langle \fl \rangle \Nm(\fl) \mathbf{T}_{\fl^{n-1}}
\]
can be used to
 determine the image of the Hecke operators attached to powers of $\fl$.
A straightforward inductive argument now shows that
\[
\begin{array}{lll}
\tilde \pi_f(\mathbf{T}_{\fl^n}) &=& 
\displaystyle \sum_{j=0}^n  \psi(\fl)^j \left(1 + \varepsilon \log_p \Nm(\fl) \left( j\cdot \lambda    + (n-j)\cdot \mu   \right)  \right) \\
& =& 
\displaystyle \sum_{I \mid \fl^n} \psi(I) \left(1 + \varepsilon \left( \lambda \log_p \Nm(I) + \mu \log_p \Nm\left(\frac{\fl^n}{I}\right) \!\! \right) \!\! \right) 
\end{array}
\]

\item For the nearly ordinary Hecke operators at $p$, it follows from \cite[Proposition 2.3]{Hid89} that 
\[
\tilde \pi_f(\mathbf{T}(p^n)) =\tilde  \pi_f(\mathbf{T}(p))^n = 1.
\]
\end{itemize}

\par We are now ready to compute the Fourier coefficients of the anti-parallel deformation $\cF$. Recall that the algebraic notion of $q$-expansions (cf.~Hida \cite[Chapter 4]{Hid04}) gives a tuple of power series $\cF_i(q)$ indexed by a set $\ft_i$ of integral ideals representing the classes in the narrow ideal class group of $F$:
\begin{equation}
\cF_i(q) \ = \ a_0(\ft_i)\ + \sum_{\nu \, \in \, (\ft_i)_+} a_{\nu }\, q^{\nu}, \qquad \qquad a_{\nu} \in E [\varepsilon]/(\varepsilon^2) ,
\end{equation}
where   as usual $(\ft_i)_+$ denotes the subset of totally positive elements of $\ft_i$. Abbreviate the $q$-expansion corresponding to the inverse different $\fd^{-1}$  by $\cF(q)$. 
The main result of this section is:
\begin{theorem}
\label{thm:alice}
 The anti-parallel family
 \begin{equation}
\cF(q) \ = \ \sum_{\nu \, \in \, \fd^{-1}_+} a_{\nu} \, q^{\nu}, \qquad a_{\nu} \in E [\varepsilon]/(\varepsilon^2) ,
\end{equation}
has Fourier  coefficients given, to first order, by
\begin{equation}
\label{eqn:Fourier_coefficients_F}
a_{\nu} =  \displaystyle 
\sum_{I \mid (\nu) \fd} \psi(I) \left(1 + \varepsilon \left(-\log_p(\nu) + \frac{\sL_1(\psi)}{\sbL(\psi)} \log_p \Nm(I) + \frac{\sL_2(\psi)}{\sbL(\psi)} \log_p \Nm\left( \frac{(\nu)\fd}{I}\right) \!\!\right)\!\! \right),
\end{equation}
for all $\nu$ that are relatively prime to $p$.
Furthermore $a_{p^m\nu} = a_{\nu}$ for all $m\ge 1$.
\begin{proof}
The family $\cF$ is $p$-adically cuspidal, so   the constant term vanishes. To compute the higher Fourier coefficients, we compute that, for any $\nu \in \fd_+^{-1}$ such that $p \nmid \nu$,  
\begin{eqnarray*}
\tilde \pi_f(\mathbf{T}_{(\nu)\fd}) &=& \prod_{\fl^n \mid \! \mid (\nu) \fd} \tilde \pi_f \left(\mathbf{T}_{\fl^n}\right) \\
&=&
\displaystyle 
\sum_{I \mid (\nu) \fd} \psi(I) \left(1 + \varepsilon \left(\frac{\sL_1(\psi)}{\sbL(\psi)} \log_p \Nm(I) + \frac{\sL_2(\psi)}{\sbL(\psi)} \log_p \Nm\left( \frac{(\nu) \fd}{I}\right)\!\! \right) \!\! \right)
\end{eqnarray*}
to first order. To determine the Fourier coefficients from this value,  the  $p$-adic interpolation properties of the coefficients  $a_{\nu}$, and the density of classical forms, are used
to reduce to the relations between Fourier coefficients and the Hecke algebra proved for classical forms in \cite{Hid91}. 

Consider the rigid analytic fiber of the formal scheme attached to $\cT^{\nord}$. For a sufficiently small affinoid neighbourhood $V=\Spm(A_V)$ of the point corresponding to the morphism $\pi_f$, there is a rigid analytic family $\cF_{V}=\sum_{\nu} a_{V, \nu}q^{\nu}$ with normalised Fourier coefficients in $a_{V,\nu} \in A_V$, specialising to $\cF$ in the anti-parallel direction. 
By Hida's Control Theorem, there is a Zariski-dense set of points in $V$ corresponding to systems of Hecke eigenvalues $\pi_g \colon A_V\to \bar \Q_p$
of classical modular forms $g$ of weight $(\bfk_g, \bfw_g)$ and with   Fourier coefficients $a_{g, \nu}$ given by the image of $a_{V,\nu}$ under $\pi_g$.
Combining the relations for classical forms proved in \cite[Eqn. (2.3) et seq./Eqn. (1.5)]{Hid91}, one obtains
\begin{equation*}
\begin{array}{lll}
\pi_g(\mathbf{T}_{(\nu) \fd}) &=& a_{g, \nu} \cdot \nu^{\bfv_g} \\
\pi_g(\mathbf{T}_{p^m}) &=& a_{g,p^m} 
\end{array}
\end{equation*}
where as before $\bfv_g = \bfw_g - \bfk_g + \bft$.  The quantity $\nu^{\bfv_g}$ may be identified with the weight $(\bfk_g, \bfw_g)$-specialisation of the universal character $\kappa^{\univ}$ evaluated at $(\nu, 1)$. The image of $\kappa^\univ (\nu, 1)$ under the morphism \eqref{eqn:antiparallel direction} defining  the anti-parallel direction  is given by 
\[
1 + \varepsilon \log_p(\nu) 
\]
so that the density of classical points in $V$ implies that 
\[
a_{\nu} = (1 - \varepsilon \log_p(\nu)) \tilde \pi_f(\mathbf{T}_{(\nu) \fd}). 
\] 
The result follows.
\end{proof}
\end{theorem}

\section{Diagonal restrictions and RM values}
\label{sec:generating-series}

This section  explains how to parlay Theorem 
\ref{thm:alice}
of \S~\ref{sec:deformation} into a proof of Theorem C. In a nutshell, the generating series of Theorem C is obtained from the ordinary projection of the diagonal restriction of a modification of the anti-parallel cuspidal deformation $\cF$ described in Theorem \ref{thm:alice}.

\par Retain the setup of \S~\ref{subsec:q-exp}. Namely, $F$ denotes a real quadratic field, and $\psi$ is a totally odd unramified character. Let $D$ be the discriminant of $F$, with ring of integers $\cO_F$, set of integral ideals $\mathscr{I}_F$, and different ideal $\fd$. The notation $\fd^{-1}_+$ is used for the subset of totally positive elements of the inverse different $\fd^{-1}$. Write $\Nm$ and $\Tr$ for the norm and trace functions from $F$ to $\Q$.  If $\tau\in \cH_p^D$ is an RM point of discriminant $D$,  denote by $\fa_\tau\in {\rm Cl}(D)$ the narrow ideal class attached to $\tau$. If $J$ is a rigid cocycle, then
\[
J[\psi]:= \prod_{\tau\in \SL_2(\Z)\backslash \cH_p^{\circ,D}}   J[\tau]^{\psi(\fa_\tau)} \in \C_p^\times \otimes \Q(\psi).
\]

\begin{remark}
\label{rem:ramified}
The character $\psi$ was assumed to be  unramified for simplicity, and it would be interesting to generalise the arguments to the case of an arbitrary totally odd ring class character
\begin{equation}
 \psi : \Cl(D) \ \lra \ \C_p^{\times}
\end{equation}
of discriminant $D = f^2 D_0$ with $D_0$ fundamental, and $(p,f)=1$. The deformations studied in \S~\ref{sec:deformation} are not sensitive in an essential way to this additional ramification. Moreover, a version of Lemma \ref{lemma:bijections} for non-trivial conductors can be found in \cite{LV19}, and the explicit formula \eqref{eqn:TnJ-formula} continues to hold. One may therefore expect that Proposition \ref{thm:ord-proj} is amenable to this generalisation via the strategy of this paper, provided that the left hand side of the equality is replaced by the series obtained by taking the trace to level $\Gamma_0(p)$:
\begin{equation}
\label{eqn:ramified}
 \mathrm{Tr}^{\Gamma_0(fp)}_{\Gamma_0(p)}\left( e^{\rm ord}(\partial_{\varepsilon}f^+) \right) \quad \in M_2(\Gamma_0(p)).
\end{equation}
\end{remark}

\subsection{The RM values of the winding cocycle}
\label{subsec:RMvalues}

In contrast the approach of \cite{darmon-dasgupta}, the calculations below build on the viewpoint of rigid (theta) cocycles introduced in \cite[\S 3]{darmon-vonk-borcherds}, by making essential use of the  winding cocycle $J_w$ of the prequel \cite{DPV1}, some of whose properties were already recalled in \S~\ref{sec:winding-cocycle}. This section describes  some further results from \cite{DPV1} concerning its RM values. The  first key result is an explicit formula for $T_nJ_w[\tau]$, which was established in \cite[Theorem 2.9]{DPV1}.

\par In order to state it, choose, for any  integer $n\geq 1$ and any RM point $\tau$ in $\cH_p$, a finite set $M_n(\tau)$ of representatives  for the double coset space $\SL_2(\Z) \backslash M_2(\Z)_n / \Gamma_\tau$, where
\[
M_2(\Z)_n  := \{\alpha\in M_2(\Z) \mbox{ with } \det(\alpha) = n\}, \qquad \quad \Gamma_\tau = \Stab_{\SL_2(\Z)}(\tau).
\]
In other words
\begin{equation}
M_2(\Z)_n = \bigsqcup_{\delta  \in M_n(\tau)} \SL_2(\Z) \cdot \delta  \cdot \Gamma_\tau.
\end{equation}
 Let $\tilde\Gamma:= \GL_2^+(\Z[1/p])$ be the group of invertible matrices over $\Z[1/p]$ with positive determinant.
\begin{theorem}
\label{thm:dpv1-Jw}
 Let  $n \geq 1$ be an integer coprime to $p$. Then 
 \begin{equation}
\label{eqn:TnJ-formula}
T_nJ_w[\tau] = 
\prod_{\delta  \in M_n(\tau)}  
\prod_{\substack{w \,\in\, \tilde{\Gamma} \delta  \tau \\ v_p(w) \,=\, 0}} w^{[0,\infty]\cdot (w',w)}.
\end{equation}
\begin{proof}
See \cite[Theorem 2.9]{DPV1}.
\end{proof}
\end{theorem}

Lemma \ref{lemma:bijections} below recalls the existence of a bijection between ``level $n$'' sets of RM points and ideals that was constructed in \cite[\S~1]{DPV1}. Define the (multi)set
\begin{equation}
\RM_n^{+}(\tau) := \bigsqcup_{\delta  \in M_n(\tau)} \left\{ w \in \tilde\Gamma \delta \tau \ : \  
\begin{array}{c}
 w > 0 > w' \\
 v_p(w) = 0, \, v_p(\mathrm{disc}(w)) \leq v_p(n) 
\end{array}
\right\}
\end{equation}
where $\mathrm{disc}(w)$ is defined to be the discriminant of a primitive integral quadratic form that has $w$ as a root. Similarly, define $\RM_n^-(\tau)$ as above, with the condition $w > 0 > w'$ replaced by $w' > 0 > w$. 

\begin{remark}
Note that an RM point $w$ may appear several times in the set $\RM_n^+(\tau)$, and the multiplicity with which it does is a subtle actor in the bijections discussed below. It is therefore important to use a disjoint union in this definition. The nature of the matrices $\delta$, which index the multiplicity with which an RM point $w$ arises, was made clearer in the proof of \cite[Lemma 1.9]{DPV1}.
\end{remark}

The sets $\RM_n^{\pm}(\tau)$ play a crucial role in the explicit formulae for the Fourier coefficients of the diagonal restrictions of the Eisenstein family $\cE$ investigated in \cite{DPV1}. It will be observed below that they appear again in the analysis of the diagonal restriction of the anti-parallel family $\cF$ studied in \S~\ref{sec:deformation}. 
\begin{lemma}
\label{lemma:bijections}
There exist
 two bijections
\[
\begin{array}{ccll}
 \phi_1 \ : & \RM_n^{-}(-\tau) & \lra & \RM_n^{+}(\tau) \\ [4pt]
 \phi_2 \ : & \left\{ (I,\nu) : \begin{array}{cc}
\nu \in \fd^{-1}_+, & p \nmid I \, | \, (\nu)  \fd \\
\Tr(\nu) = n, & I \sim (1,\tau)
\end{array} \right\} & \lra & \RM_n^{+}(\tau) \\
\end{array}
\]
such that, after writing $\nu = p^m \nu_0$, we have
\[
\begin{array}{lll}
\phi_1(w) &=& -w, \\
\phi_2(I,\nu) &=& \nu_0\sqrt{\Delta}/\Nm(I).
\end{array}
\]
\begin{proof}
A bijection $\phi_1$ as required may be constructed by letting $W_{\infty}$ be a diagonal matrix with eigenvalues $1$ and $-1$. If $w \in \tilde\Gamma\delta\tau$, then  
\[
-w = W_{\infty}w \ \in \ \tilde\Gamma\delta'(-\tau)
\]
where $\delta' \in M_n(\tau)$ is the double coset representative of $W_{\infty}\delta W_{\infty}$. To obtain a bijection $\phi_2$ as above, one first uses a bijection
\[
\Phi \ : \ \ 
\left\{ (I,\nu) : \begin{array}{cc}
\nu \in \fd^{-1}_+, & I \, | \, (\nu) \fd \\
\Tr(\nu) = n, & I \sim (1,\tau)
\end{array} \right\} \ \ \lra \ 
\bigsqcup_{\delta \in M_n(\tau)} \left\{ w \in \SL_2(\Z) \delta\tau \ : \  
 w > 0 > w' 
\right\}
\]
which satisfies $\Phi(I,\nu) = \nu\sqrt{\Delta}/\Nm(I)$. Such a bijection was constructed in \cite[Lemma 1.9]{DPV1}. Note that the source of $\Phi$ is almost equal to the source of the desired bijection, minus the condition $p \nmid I$. Under the bijection $\Phi$, the condition $p \nmid I$ is equivalent to the condition that $w = w_0p^m$ for some $m\geq 0$ and $p\nmid w_0$. The map $w \mapsto w_0$ then defines a bijection between
\[
\left\{ w \in \SL_2(\Z) \delta\tau \ : \  
\begin{array}{c}
 w > 0 > w' \\
 w = w_0p^m, \, p \nmid w_0, \, m \geq 0
\end{array}
\right\} ,
\]
and the set
\[
\left\{ w \in \tilde\Gamma \delta\tau \ : \  
\begin{array}{c}
 w > 0 > w' \\
 v_p(w) = 0, \, v_p(\mathrm{disc}(w)) \leq v_p(n) 
\end{array}
\right\},
\]
so that the result follows by definition of $\RM_n^+(\tau)$. 
\end{proof}
\end{lemma}

\subsection{Derivatives of diagonal restrictions. }
The modular generating series for the RM values of the winding cocycle that is the subject of Theorem C will be constructed from three different analytic families that specialise to the Eisenstein series of parallel weight one. More specifically, the anti-parallel cuspidal family from \S~\ref{subsec:q-exp}, and the two Eisenstein families of Lemma \ref{lemma:eisenstein-family}:
\[
E^{(p)}_{1+ \varepsilon}(1,\psi) \quad \mbox{and } \quad E^{(p)}_{1+ \varepsilon}(\psi,1)
\]

\par The \textit{modularity} of the generating series of Theorem C will follow from two simple results: 
\begin{enumerate}
\item The vanishing of the diagonal restriction of $E^{(p)}_1(1,\psi)$, the $p$-stabilisation of the parallel weight one Eisenstein series (for which the shorthand $f$ was used in \S~\ref{sec:deformation}), 
\item For any analytic family of $p$-adic modular forms whose specialisation vanishes, the specialisation of its \textit{derivative} is also a $p$-adic modular form. 
\end{enumerate} 
These results were also used in \cite{DPV1}, where full proofs may be found. Since they play an important role in the argument, they will be briefly reviewed here.

\begin{lemma}
\label{lemma:diagonal_res_vanishing}
Suppose $p$ is inert in $F$, and $\psi$ is an odd unramified character of $F$. Then 
\[
E_{1}^{(p)}(1,\psi)(z,z) = 0.
\] 
\begin{proof}
Recall that the diagonal restriction of any Hilbert modular form with Fourier coefficients $a_{\nu}$ has the following $q$-expansion:
\begin{equation}
\label{eqn:diag-res-general}
a_0 + \sum_{n \geq 1} \sum_{\substack{\nu \in \fd_+^{-1} \\ \Tr(\nu) \, = \, n}} a_{\nu}  \ q^n.
\end{equation}
For the Eisenstein series $E_1^{(p)}(1,\psi)$, the Fourier coefficient $a_{\nu}$ is equal to
\[
4 \sum_{p \nmid I \mid (\nu)\fd} \psi(I)
\]
For any ideal $I$ in the index set of this summation, we may write $IJ(p^e) = (\nu)\fd$ for some uniquely determined ideal $J$ coprime to $p$, since $p$ is inert in $F$. The conjugate $J'$ then defines an ideal coprime to $p$, dividing $(\nu')\fd$. Observe that, since $\psi$ is odd, we have
\[
\psi(J') = \psi(J)^{-1} = \psi(I)\psi(\fd)^{-1} = - \psi(I),
\] 
and it  follows that $a_{\nu} = - a_{\nu'}$. Therefore the diagonal restriction must vanish. 
\end{proof}  
\end{lemma}
 
 \par The three analytic families that specialise to $E_1^{(p)}(1,\psi)$ therefore give families of diagonal restrictions that specialise to zero. It is easy to see that the specialisation of the \textit{derivative} of each of these families of diagonal restrictions is a $p$-adic modular form of weight two. The following result, appearing as Lemma 2.1 in \cite{DPV1}, ascertains that it is even overconvergent, though this is not used in what follows.

\begin{lemma}
\label{lemma:diagonal_res_derivative}
Suppose $\cF(t)$ is a family of overconvergent forms of weight $\kappa(t)$, indexed by a parameter $\, t$ on a closed rigid analytic disk $D$ in weight space. Suppose that
\begin{itemize}
\item the disk $D$ is centred at an integer $k = \kappa(0) \in \Z$,
\item the specialisation vanishes $\cF(0)=0$.
\end{itemize}
Then the derivative $\partial_t\mathcal{F}(0)$ is an overconvergent modular form of weight $k$.
\end{lemma}

\subsection{Proof of Theorem C }
Theorem \ref{thm:alice} will now be used to construct the modular generating series $G_\tau$ of Theorem C, and calculate its constant term. The argument involves three main steps:
\begin{enumerate}
\item The definition of the power series $\partial\cF^+_{\psi}$, a combination of the $q$-expansions of the first derivatives of the anti-parallel cuspidal family $\cF_{\psi}$ of Theorem \ref{thm:alice} and a parallel Eisenstein family $\cE_{\psi}$;
\item The computation of its diagonal restriction $\partial f^+_{\psi}$;
\item The computation of its ordinary projection $e^{\ord}(\partial f^+_{\psi})$. 
\end{enumerate}
The forms constructed in these three steps lie in increasingly structured spaces: $\partial f^+_{\psi}$ is a $p$-adic modular form of weight two and tame level one, and $e^{\ord}(\partial f^+_{\psi})$ is a classical modular form on $\Gamma_0(p)$. The power series $\partial\cF^+_{\psi}$ however lacks the modularity properties of a traditional (classical or $p$-adic) Hilbert modular form, and is perhaps best envisaged as an instance of a ``$p$-adic mock modular form", of the kind that make an appearance in \cite{DT08, darmon-lauder-rotger1} for instance.

\par The series $\cF^+_{\psi}$  is a combination of first order families of modular forms passing through the same Eisenstein series of parallel weight one in \textit{different} weight directions. Its definition was dictated by the algebraic shape of the Fourier coefficients of the anti-parallel family $\cF_{\psi}$ arising from Theorem \ref{thm:alice}, as it causes the desired algebraic cancellation. Precisely, define
\[
\cF^+_{\psi} \ := \ \cF_{\psi} + \cE_{\psi} \ = \ a_{0}(\cF^+_{\psi})  \ + \sum a_{\nu}(\cF^+_{\psi})  \, q^{\nu}, \qquad \qquad a_{\nu}(\cF^+_{\psi}) \in E [\varepsilon]/(\varepsilon^2) ,
\]
where the first term $\cF_{\psi}$ is the anti-parallel weight cuspidal deformation of Theorem \ref{thm:alice}. The second term $\cE_{\psi}$ is the following explicit combination of parallel weight Eisenstein families 
\begin{equation}
\cE_{\psi} := \frac{\sL_2(\psi)}{\sbL(\psi)} \left(  E_{1+\varepsilon}^{(p)}(1,\psi) \ - \ \, E_{1+\varepsilon}^{(p)}(\psi,1)\right).
\end{equation}
Recall the Gross--Stark unit $u_\psi$ attached to the odd character $\psi$, defined in \eqref{eqn:def-u-phi}. Henceforth, the unit $u_\psi$ is identified with its image under the $p$-adic embedding $\tilde\alpha_2$ in order to lighten the notations and view it as an element of $F_p^\times \otimes \Q(\psi)$, to which the $p$-adic logarithm $\log_p$ may be unambiguously applied.
\begin{proposition}
The Fourier coefficients of $\cF^+_{\psi}$ are given by
$$
\begin{array}{lllll}
a_0(\cF^+_{\psi}) &=& \frac{\varepsilon}{2} \cdot  \log_p(u_\psi).
\\ [3pt]
a_{\nu}(\cF^{+}_{\psi})
&=& \displaystyle 
\sum_{I \mid (\nu_0) \fd} \psi(I) \left( 1-\varepsilon\log_p \left( \frac{\nu_0}{\Nm(I)} \right)\!\! \right) & &
\end{array} 
$$
\begin{proof}
Since $\psi$ is odd, we have $\psi(\fd)=-1$. The constant term of $\cE_{\psi}$, given by \eqref{eqn:Fourier_E_eta_psi}, is therefore
\[
a_0(\cE_{\psi}) =   \frac{\sL_2(\psi)}{\sbL(\psi)} \left( a_0(1,\psi)  -  a_0(\psi,1)  \right) 
= \frac{\varepsilon}{4}\cdot  \frac{\sL_2(\psi)}{\sbL(\psi)} (L_p'(\psi,0) + L_p'(\psi^{-1},0)).
\]
By the Gross--Stark theorem \eqref{eqn:Gross--Stark} 
\begin{equation}
L_p'(\psi,0) = \sbL(\psi) L(\psi,0) = \sbL(\psi^{-1}) L(\psi^{-1},0) = L_p'(\psi^{-1},0),
\end{equation}
and hence, using Lemma \ref{lemma:from-L-to-u}, we obtain
\[
a_0(\cE_{\psi}) 
=  \frac{\varepsilon}{2}\cdot  \frac{\sL_2(\psi)}{\sbL(\psi)} L_p'(\psi,0) 
=  \frac{\varepsilon}{2}\cdot  \sL_2(\psi) L(\psi,0) = \frac{\varepsilon}{2} \cdot \log_p(u_\psi).
\]

\par At $\nu\ne 0$, the Fourier coefficient of $\cE_{\psi}$ is given by
$$ a_\nu(\cE_{\psi}) = 
 \displaystyle \frac{\sL_2(\psi)}{\sbL(\psi)} 
\sum_{I \mid (\nu_0) \fd} \psi(I) \varepsilon \left(\log_p(\Nm(I)   -\log_p \left( \frac{\nu_0}{\Nm(I)} \right)\!\! \right).$$
Combining this with the  formula for  the Fourier coefficients of the anti-parallel deformation $\cF_{\psi}$ given in Theorem \ref{thm:alice}, gives the required identity
\begin{equation}
\label{eqn:Fourier_F+}
 a_{\nu}(\cF^{+}_{\psi}) 
\quad = \quad \displaystyle 
\sum_{I \mid (\nu_0) \fd} \psi(I) \left( 1-\varepsilon\log_p \left( \frac{\nu_0}{\Nm(I)} \right)\!\! \right)
\end{equation}
 \end{proof}
 \end{proposition}
 
\par Next, we consider the \textit{diagonal restriction} $f^+_{\psi}$ of the series $\cF^+_{\psi}$, defined as the sum of the diagonal restrictions of the families $\cF_{\psi}$ and $\cE_{\psi}$. Its derivative with respect to $\varepsilon$ is modular. More specifically:
\begin{proposition}
The power series 
\begin{equation}
\label{eqn:q-exp-G}
 \partial f^+_{\psi}(q) \ = \ \frac{1}{2} \log_p(u_\psi) \ - \ \sum_{n \geq 1} \sum_{\substack{\nu \in \fd_+^{-1} \\ \Tr(\nu) \, = \, n}}\sum_{I \mid (\nu_0) \fd} \psi(I) \log_p \left( \frac{\nu_0\sqrt{D}}{\Nm(I)} \right) q^n. 
\end{equation}
is the $q$-expansion of a $p$-adic modular form of weight two and tame level one.
\begin{proof}
Lemma \ref{lemma:diagonal_res_vanishing} implies that the diagonal restriction $f^+_{\psi}$ vanishes at $\varepsilon = 0$, so that the derivative $\partial f^+_{\psi}$ is a $p$-adic modular form (by Lemma \ref{lemma:diagonal_res_derivative} it is even overconvergent). The statement about its $q$-expansion follows by \eqref{eqn:diag-res-general} from the observation that $\partial f^+_{\psi}(q)$ differs from the desired result by
\[
\sum_{\substack{\nu \in \fd_+^{-1} \\ \Tr(\nu) \, = \, n}} \psi(I) \log_p(\sqrt{D}),
\] 
which is proportional to the $n$-th Fourier coefficient of the diagonal restriction of the Hilbert Eisenstein series $E_1^{(p)}(1, \psi)$, and is therefore identically zero by Lemma \ref{lemma:diagonal_res_vanishing}. 
\end{proof}
\end{proposition}

\par Finally, we explicitly compute the \textit{ordinary projection} of the $p$-adic modular form $\partial f^+_{\psi}$. This ordinary projection is a classical modular form in $M_2(\Gamma_0(p))$, and its Fourier coefficients can be related to the RM values of the winding cocycle $J_w$, using the explicit formula for the latter stated in \S~\ref{subsec:RMvalues}.

\begin{proposition}
\label{thm:ord-proj}
The ordinary projection of the $p$-adic modular form $\partial f^+_{\psi}$ is a classical modular form in the space $M_2(\Gamma_0(p))$. Its $q$-expansion is given by:
\begin{equation}
2e^{\rm ord}(\partial f^+_{\psi}) \,=\, \log_p(u_\psi)
 \ - \ \sum_{n \geq 1} \log_p \left( T_n J_w[\psi] \right) q^n.
\end{equation}
\begin{proof}
Note that the ordinary projection is classical of level $\Gamma_0(p)$, by Coleman's classicality theorem. The statement about the constant term follows from \eqref{eqn:Fourier_F+}. For any $n \geq 1$, the bijection $\phi_2$ of Lemma \ref{lemma:bijections} allows us to rewrite the $n$-th Fourier coefficient of $2\partial f^+_{\psi}$ appearing in \eqref{eqn:q-exp-G} in terms of the level $n$ sets of RM points $\RM_n^{\pm}(\tau)$. Since $\psi(\tau) = - \psi(-\tau)$, this Fourier coefficient is given by
\begin{eqnarray*}
2a_n &=& 
\sum_{\tau\in \SL_2(\Z)\backslash \cH_p^{\circ,D}} \psi(\tau)
\left( 
\sum_{w \,\in\, \RM_n^+(\tau)} \log_p(w) \ - \sum_{w \,\in\, \RM_n^+(-\tau)} \log_p(w) \right) \\
&=& 
\sum_{\tau\in \SL_2(\Z)\backslash \cH_p^{\circ,D}} \psi(\tau)
\left( 
\sum_{w \,\in\, \RM_n^+(\tau)} \log_p(w) \ - \ \sum_{w \,\in\, \RM_n^-(\tau)} \ \log_p(w) 
\right)
\end{eqnarray*}
where the second equality follows from the existence of a bijection $\phi_1$ as in Lemma \ref{lemma:bijections}. Let $n \geq 1$ coprime to $p$, then the $n$-th coefficient of the ordinary projection of $2\partial f^+_{\psi}$ is given by
\begin{eqnarray*}
2a_n^{\rm ord} = 2 \lim_{m \to \infty} a_{np^{2m}} 
 &=& \!\! 
 \sum_{\tau\in \SL_2(\Z)\backslash \cH_p^{\circ,D}} \psi(\tau)
 \sum_{\delta  \in  M_n(\tau)} 
 \sum_{\substack{w \, \in \, \tilde{\Gamma} \delta \tau \\ v_p(w) \, = \, 0}} 
 \bigl([0,\infty]\cdot (w',w) \bigr)\log_p(w) \\
 &=& \ \log_p (T_nJ_w[\psi])
\end{eqnarray*}
where the last equality uses the explicit formula \eqref{eqn:TnJ-formula}. 

\par To obtain the statement for all $n \geq 1$, note that 
\[
T_n \ \longmapsto \ \log_p (T_nJ_w[\psi]),
\]
is a linear function from the weight two Hecke algebra of level $\Gamma_0(p)$, since the Hecke action on $\rH^1(\Gamma,\cA^{\times}/\C_p^{\times})$ factors through it, and hence there exists $f \in M_2(\Gamma_0(p))$ with higher Fourier coefficients as in the statement. By what we showed, the $n$-th Fourier coefficients of $f$ and $2e^{\rm ord}(\partial f^+_{\psi})$ agree when $n$ is coprime to $p$. Their difference must therefore be an oldform, and hence zero. 
\end{proof}
\end{proposition}

We are now ready to prove Theorem C of the introduction:

\begin{theorem} 
Let $D$ be a fundamental discriminant and let 
$\tau\in \cH_p$ be an  RM point of discriminant $D$.
There is a   classical modular form $G_\tau$
of weight two on $\Gamma_0(p)$ with $p$-adic Fourier coefficients,  whose $q$-expansion is given by
$$ G_\tau(q) =  \log(u_\tau) + \sum_{n=1}^\infty \log((T_n J_w)[\tau]) q^n,$$
where $ \log \ : \ \cO_{\C_p}^\times \lra \C_p $
 is  the $p$-adic logarithm.  
 The modular form $G_\tau$  is non-trivial if and only if $\Q(\sqrt{D})$ does not admit a unit of 
norm $-1$.
\begin{proof}
Let $H$ be the narrow class field of $\Q(\sqrt{D})$.
Proposition \ref{thm:ord-proj} produces, for each odd character $\psi$ of ${\rm Gal}(H/F)$, a classical modular form in $M_2(\Gamma_0(p))$ with $q$-expansion  in $C_p[[q]]$ given by
$$ G_\psi(q) =   \log_p(u_\psi)
 \ - \ \sum_{n \geq 1} \log_p \left( T_n J_w[\psi] \right) q^n.
$$
The assignment $\psi \mapsto G_\psi(q)$ extends by linearity to a map on the linear span of the odd characters, which is the space of odd functions on ${\rm Gal}(H/F)$. Let $\psi$ be the odd indicator function on the class of $\tau$, which is equal to  $1$ on $[\tau]$, to $-1$ on $[-\tau] = [\sigma_\infty \tau]$,  where $\sigma_\infty\in {\rm Gal}(H/F)$ is complex conjugation,
and vanishes on all the other ${\rm Pic}^+(\cO_F)$-translates of $\tau \in \SL_2(\Z)\backslash \cH_p^D$. 
With this choice of $\psi$, we have
$$ 
\log_p(u_\psi) = \log_p(u_\tau) - \log_p(\sigma_\infty u_\tau) = 2 \log_p(u_\tau),  $$
and
$$
 \log_p(T_n J_w[\psi])  = \log_p(T_n J_w [\tau]) - \log_p(T_n J_w[-\tau]) =  2 \log_p(T_n J_w[\tau]).$$
The modular form $G_\tau$ of Theorem C
is obtained by setting
 $$ G_\tau = \frac{1}{2} G_\psi.$$
\end{proof}
\end{theorem}

 

\end{document}